\documentclass[11pt,a4paper,fleqn]{article}
\usepackage{a4wide,amsfonts,amsmath,latexsym,amssymb,euscript,graphicx,units,mathrsfs}

\usepackage{graphicx}
\usepackage{color}
\usepackage{amssymb}
\usepackage{amssymb}
\usepackage[T1]{fontenc}
\usepackage{latexsym}
\usepackage{xypic}
\usepackage{eufrak}
\usepackage{euscript}
\usepackage{amsfonts,amsmath}
\usepackage{verbatim}
\usepackage{fancyhdr}
\usepackage[english]{babel}
\usepackage{mathrsfs}
\usepackage{units}
\usepackage{stmaryrd}

\newtheorem{prop}{Proposition}[section]
\newtheorem{cor}[prop]{Corollary}
\newtheorem{lemme}[prop]{Lemma}
\newtheorem{rem}[prop]{Remark}

\newtheorem{thm}[prop]{Theorem}

\renewcommand{\geq}{\geqslant}
\def\leq{\leqslant}

\newcommand{\N}{\mathbb{N}}
\newcommand{\Z}{\mathbb{Z}}
\newcommand{\R}{\mathbb{R}}

\def\HH{\EuFrak H}

\def\lcr{\left[}

\def\1{{\mathbf{1}}}

\def\rcr{\right]}

\def\sk{{\mathbb{D}}}

\def\1{{\mathbf{1}}}
\def\0.5{{\frac{1}{2}}}

\def\rcr{\right]}

\newcommand{\fin}
{ \vspace{-0.6cm}
\begin{flushright}
\mbox{$\Box$}
\end{flushright}
\noindent }

\newcommand{\qed}{\nopagebreak\hspace*{\fill}
{\vrule width6pt height6ptdepth0pt}\par}

\begin{document}

\begin{center}
{\Large{\bf Stein's method on Wiener chaos}}\\~\\
by Ivan Nourdin\footnote{Laboratoire de Probabilit{\'e}s et
Mod{\`e}les Al{\'e}atoires, Universit{\'e} Pierre et Marie Curie,
Bo{\^\i}te courrier 188, 4 Place Jussieu, 75252 Paris Cedex 5,
France, {\tt ivan.nourdin@upmc.fr}} and Giovanni
Peccati\footnote{ Laboratoire de Statistique Th\'eorique et
Appliqu\'ee, Universit\'e Pierre et Marie Curie, 8\`eme \'etage,
b\^atiment A, 175 rue du Chevaleret, 75013 Paris, France,
{\tt giovanni.peccati@gmail.com}}\\
{\it University of Paris VI}\\~\\
\small{{\it Revised version}: May 10, 2008}
\end{center}

{\small \noindent {\bf Abstract:} We combine Malliavin calculus
with Stein's method, in order to derive explicit bounds in the
Gaussian and Gamma approximations of random variables in a fixed
Wiener chaos of a general Gaussian process. Our approach
generalizes, refines and unifies the central and non-central limit
theorems for multiple Wiener-It\^{o} integrals recently proved (in
several papers, from 2005 to 2007) by Nourdin, Nualart,
Ortiz-Latorre, Peccati and Tudor. We apply our techniques to prove
Berry-Ess\'een bounds in the Breuer-Major CLT for subordinated
functionals of fractional Brownian motion. By using the well-known
Mehler's formula for Ornstein-Uhlenbeck semigroups, we also
recover a technical result recently proved by Chatterjee,
concerning the Gaussian approximation of functionals of
finite-dimensional Gaussian vectors. }
\\


\\

\section{Introduction and overview}
\subsection{Motivations}
Let $Z$ be a random variable whose law is
absolutely continuous with respect to the Lebesgue measure (for
instance, $Z$ is a
standard Gaussian random variable 
or a Gamma random variable). Suppose that $\{Z_n: n \geq 1\}$ is a
sequence of random variables converging in distribution towards
$Z$, that is:
\begin{equation}\label{naive}
\mbox{for all $z\in\R$,}\quad P(Z_n\leq z)\longrightarrow P(Z\leq
z)\quad\mbox{as $n\rightarrow\infty$.}
\end{equation}
It is sometimes possible to associate an explicit uniform bound
with the convergence (\ref{naive}), providing a global description
of the error one makes when replacing $P(Z_n\leq z)$ with $P(Z\leq
z)$ for a fixed $n\geq 1$. One of the most celebrated results in
this direction is the following \textsl{Berry-Ess\'een Theorem}
(see e.g. Feller \cite{Feller book} for a proof), that we record
here for future reference:

\begin{thm}[Berry-Ess\'een]\label{T : BE}
Let $(U_j)_{j\geq 1}$ be a sequence of independent and identically
distributed random variables, such that $E(|U_j|^3)=\rho<\infty$,
$E(U_j)=0$ and $E(U_j^2)=\sigma^2$. Then, by setting $Z_n$ $=$
$\frac{1}{\sigma\sqrt{n}}\sum_{j=1}^n U_j,\, n\geq 1,$ one has
that $Z_n \,\,\,{\stackrel{{\rm Law}}{\longrightarrow}}\,\,\,
Z\sim\mathscr{N}(0,1)$, as $n\rightarrow\infty$, and moreover:
\begin{equation}\label{BerriEss\'een}
\sup_{z\in\R}|P(Z_n\leq z)-P(Z\leq
z)|\leq\frac{3\rho}{\sigma^3\sqrt{n}}.
\end{equation}
\end{thm}

The aim of this paper is to show that one can combine
\textsl{Malliavin calculus} (see e.g. \cite{Nbook}) and
\textsl{Stein's method} (see e.g. \cite{Chen_Shao_sur}), in order
to obtain bounds analogous to (\ref{BerriEss\'een}), whenever the
random variables $Z_n$ in (\ref{naive}) can be represented as
functionals of a given Gaussian field. Our results are general, in
the sense that (\textbf{i}) they do not rely on any specific
assumption on the underlying Gaussian field, (\textbf{ii}) they do
not require that the variables $Z_n$ have the specific form of
partial sums, and (\textbf{iii}) they allow to deal (at least in
the case of Gaussian approximations) with several different
notions of \textsl{distance} between probability measures. As
suggested by the title, a prominent role will be played by random
variables belonging to a \textsl{Wiener chaos} of order $q$
($q\geq 2$), that is, random variables having the form of a
\textsl{multiple stochastic Wiener-Itô integral} of order $q$ (see
Section \ref{SS : Elements Mall} below for precise definitions).
It will be shown that our results provide substantial refinements
of the central and non-central limit theorems for multiple
stochastic integrals, recently proved in \cite{NouPec07} and
\cite{NP}. Among other applications and examples, we will provide
explicit Berry-Ess\'een bounds in the Breuer-Major CLT (see
\cite{BM}) for fields subordinated to a fractional Brownian
motion.

Concerning point (\textbf{iii}), we shall note that, as a
by-product of the flexibility of Stein's method, we will indeed
establish bounds for Gaussian approximations related to a number
of distances of the type
\begin{equation}\label{distances}
d_\mathscr{H}(X,Y)=\sup\left\{ \left|E(h(X))-E(h(Y))\right|\;:\;
h\in\mathscr{H}\right\},
\end{equation}
where $\mathscr{H}$ is some suitable class of functions. For
instance: by taking $\mathscr{H}=\{h : \Vert h\Vert _{L}\leq 1 \},
$ where $\left\Vert \cdot \right\Vert _{L}$ denotes the usual Lipschitz
seminorm, one obtains the \textsl{Wasserstein} (or
\textsl{Kantorovich-Wasserstein}) distance; by taking
$\mathscr{H}=\{h : \Vert h\Vert _{BL}\leq 1 \}, $ where
$\left\Vert \cdot \right\Vert _{BL}$ $=$ $\left\Vert \cdot
\right\Vert _{L}+\left\Vert \cdot \right\Vert _{\infty}$, one
obtains the \textsl{Fortet-Mourier} (or \textsl{bounded
Wasserstein}) distance; by taking $\mathscr{H}$ equal to the
collection of all indicators $\mathbf{1}_B$ of Borel sets, one
obtains the \textsl{total variation distance}; by taking
$\mathscr{H}$ equal to the class of all indicators functions
$\mathbf{1}_{(-\infty,z]}$, $z\in\R$, one has the
\textsl{Kolmogorov} distance, which is the one taken into account
in the Berry-Ess\'een bound (\ref{BerriEss\'een}). In what follows, we
shall sometimes denote by $d_{\rm W}(.,.)$, $d_{\rm FM}(.,.)$,
$d_{\rm TV}(.,.)$ and $d_{\rm Kol}(.,.)$, respectively, the
Wasserstein, Fortet-Mourier, total variation and Kolmogorov
distances. Observe that $d_{\rm W}(.,.)\geq d_{\rm FM}(.,.)$ and
$d_{\rm TV}(.,.)\geq d_{\rm Kol}(.,.)$. Also, the topologies
induced by $d_{\rm W}$, $d_{\rm TV}$ and $d_{\rm Kol}$ are
stronger than the topology of convergence in distribution, while
one can show that $d_{\rm FM}$ metrizes the
convergence in distribution (see e.g. \cite[Ch. 11]{Dudley book}
for these and further results involving distances on spaces of
probability measures).

\subsection{Stein's method} \label{SS : Stein Intro}

We shall now give a short account of Stein's method, which is
basically a set of techniques allowing to evaluate distances of
the type (\ref{distances}) by means of differential operators.
This theory has been initiated by Stein in the path-breaking paper
\cite{Stein_orig}, and then further developed in the monograph
\cite{Stein_book}. The reader is referred to \cite{Chen_Shao_sur},
\cite{Reinert_sur} and \cite{RinRot} for detailed surveys of
recent results and applications. The paper by Chatterjee
\cite{Chatterjee_AOP} provides further insights into the existing
literature. In what follows, we will apply Stein's method to two
types of approximations, namely Gaussian and (centered) Gamma. We
shall denote by $\mathscr{N}(0,1)$ a standard Gaussian random
variable. The centered Gamma random variables we are interested in
have the form
\begin{equation}\label{GammaCent}
F(\nu)\stackrel{\rm Law}{=}2G(\nu/2)-\nu, \quad \nu>0,
\end{equation}
where $G(\nu/2)$ has a Gamma law with parameter $\nu/2$. This
means that $G(\nu/2)$ is a (a.s. strictly positive) random
variable with density $g(x) = \frac{x^{\frac{\nu}{2}-1}{\rm
e}^{-x}}{\Gamma(\nu/2)}\mathbf{1}_{(0,\infty)}(x),$ where $\Gamma$
is the usual Gamma function. We choose this parametrization in
order to facilitate the connection with our previous paper
\cite{NouPec07} (observe in particular that, if $\nu\geq 1$ is an
integer, then $F(\nu)$ has a centered $\chi^2$ distribution with
$\nu$ degrees of freedom).

\underline{Standard Gaussian distribution.} Let $Z\sim
\mathscr{N}(0,1)$. Consider a real-valued function $h:\R\rightarrow\R$
such that the expectation $E(h(Z))$ is well-defined. The
\textsl{Stein equation} associated with $h$ and $Z$ is classically
given by
\begin{equation}\label{SteinGaussEq}
h(x)-E(h(Z)) = f'(x)-xf(x), \quad x\in\R.
\end{equation}
A solution to (\ref{SteinGaussEq}) is a function $f$ which is
Lebesgue a.e.-differentiable, and such that there exists a version
of $f'$ verifying (\ref{SteinGaussEq}) for every $x\in\R$. The
following result is basically due to Stein \cite{Stein_orig,
Stein_book}. The proof of point (i) (whose content is usually
referred as \textsl{Stein's lemma}) involves a standard use of the
Fubini theorem (see e.g. \cite{Stein64} or \cite[Lemma
2.1]{Chen_Shao_sur}). Point (ii) is proved e.g. in \cite[Lemma
2.2]{Chen_Shao_sur}; point (iii) can be obtained by combining e.g.
the arguments in \cite[p. 25]{Stein_book} and \cite[Lemma
5.1]{Chatterjee_ptrf}; a proof of point (iv) is contained in
\cite[Lemma 3, p. 25]{Stein_book}; point (v) is proved in
\cite[Lemma 4.3]{Chatterjee_AOP}.
\begin{lemme}\label{Stein_Lemma_Gauss}
\begin{enumerate}
\item[(i)] Let $W$ be a random variable. Then, $W\stackrel{\rm
Law}{=}Z\sim \mathscr{N}(0,1)$ if, and only if,
\begin{equation}\label{SteinTruc}
E[f'(W)-Wf(W)]=0,
\end{equation}
for every continuous and piecewise continuously differentiable
function $f$ verifying the relation $E|f'(Z)|$ $<$ $\infty$.
\item[(ii)] If $h(x)=\mathbf{1}_{(-\infty,z]}(x)$, $z\in\R$, then (\ref{SteinGaussEq})
admits a solution $f$ which is bounded by $\sqrt{2\pi}/4$,
piecewise continuously differentiable and such that $\|f'\|_\infty \leq 1$.
\item[(iii)] If $h$ is bounded by 1/2, then (\ref{SteinGaussEq})
admits a solution $f$ which is bounded by $\sqrt{\pi/2}$,
Lebesgue a.e. differentiable
and such that $\|f'\|_\infty
\leq 2$.
\item[(iv)] If $h$ is bounded and absolutely continuous (then,
in particular, Lebesgue-a.e. differentiable), then (\ref{SteinGaussEq}) has a
solution $f$ which is bounded and twice differentiable, and such
that $\|f\|_\infty \leq \sqrt{\pi/2} \|h-E(h(Z))\|_\infty$,
$\|f'\|_\infty \leq 2 \|h-E(h(Z))\|_\infty$ and $\|f''\|_\infty
\leq 2 \|h'\|_\infty$.
\item[(v)] If $h$ is absolutely continuous with bounded derivative,
then (\ref{SteinGaussEq}) has a solution $f$ which is twice
differentiable and such that $\|f'\|_\infty \leq \|h'\|_\infty$
and $\|f''\|_\infty \leq 2 \|h'\|_\infty$.
\end{enumerate}
\end{lemme}
We also recall the relation:
\begin{equation}\label{dollars}
2d_{\rm TV}(X,Y)=\sup\{|E(u(X))-E(u(Y))|:\|u\|_\infty\leq1\}.
\end{equation}
Note that point (ii) and (iii) (via (\ref{dollars})) imply the
following bounds on the Kolmogorov and total variation distance
between $Z$ and an arbitrary random variable $Y$:
\begin{eqnarray}\label{BoundKol}
d_{\rm Kol}(Y,Z)\leq \sup_{f\in\mathscr{F}_{\rm
Kol}}|E(f'(Y)-Yf(Y))|\\
d_{\rm TV}(Y,Z)\leq \sup_{f\in\mathscr{F}_{\rm
TV}}|E(f'(Y)-Yf(Y))| \label{boundTV}
\end{eqnarray}
where $\mathscr{F}_{\rm Kol}$ and $\mathscr{F}_{\rm TV}$ are,
respectively, the class of piecewise continuously differentiable
functions that are bounded by $\sqrt{2\pi}/4$ and such that their
derivative is bounded by 1, and the class of piecewise
continuously differentiable functions that are bounded by
$\sqrt{\pi/2}$ and such that their derivative is bounded by 2.

Analogously,
by using (iv) and (v) along with the relation
$\|h\|_L=\|h'\|_{\infty}$, one obtains
\begin{eqnarray}
d_{\rm FM}(Y,Z)\leq \sup_{f\in\mathscr{F}_{\rm
FM}}|E(f'(Y)-Yf(Y))|, \label{BoundFM}\\
d_{\rm W}(Y,Z)\leq \sup_{f\in\mathscr{F}_{\rm W}}|E(f'(Y)-Yf(Y))|,
\label{BoundWass}
\end{eqnarray}
where: $\mathscr{F}_{\rm FM}$ is the class of twice differentiable
functions that are bounded by $\sqrt{2\pi}$, whose first
derivative is bounded by 4, and whose second derivative is bounded
by 2; $\mathscr{F}_{\rm W}$ is the class of twice differentiable
functions, whose first derivative is bounded by 1 and whose second
derivative is bounded by 2.

\underline{Centered Gamma distribution.} Let $F(\nu)$ be as in
(\ref{GammaCent}). Consider a real-valued function $h:\R\rightarrow\R$
such that the expectation $E[h(F(\nu))]$ exists. The \textsl{Stein
equation} associated with $h$ and $F(\nu)$ is:
\begin{eqnarray}\label{SteinGammaEq}
h(x)-E[h(F(\nu))] = 2(x+\nu)f'(x)-xf(x),\quad x\in (-\nu,
+\infty).
\end{eqnarray}
The following statement collects some slight variations around
results proved by Stein \cite{Stein_book}, Diaconis and Zabell
\cite{DiacZab}, Luk \cite{Luk}, Schoutens \cite{Schoutens} and
Pickett \cite{These_Pickett}. It is the ``Gamma counterpart'' of
Lemma \ref{Stein_Lemma_Gauss}. The proof is detailed in Section
\ref{S : Proofs}.

\begin{lemme}\label{Stein_Lemma_Gamma}
\begin{enumerate}
\item[(i)] Let $W$ be a real-valued random variable (not necessarily
with values in $(-\nu,+\infty)$) whose law admits a density with
respect to the Lebesgue measure. Then, $W\stackrel{\rm Law}{=}
F(\nu)$ if, and only if,
\begin{equation}\label{SteinTruc2}
E[2(W+\nu)_+f'(W)-Wf(W)]=0,
\end{equation}
where $a_+:=\max(a,0)$, for every smooth function $f$ such that
the mapping $x\mapsto 2(x+\nu)_+f'(x) -xf(x)$ is bounded.
\item[(ii)] If $|h(x)|\leq c \exp(ax)$ for
every $x>-\nu$ and for some $c>0$ and $a<1/2$, and if $h$ is twice
differentiable, then (\ref{SteinGammaEq}) has a solution $f$ which
is bounded on $(-\nu,+\infty)$, differentiable and such that
$\|f\|_\infty \leq 2 \|h'\|_\infty$ and $\|f'\|_\infty \leq
\|h''\|_\infty$.

\item[(iii)] Suppose that $\nu\geq 1$ is an integer. If $|h(x)|\leq c \exp(ax)$
for every $x>-\nu$ and for some $c>0$ and $a<1/2$, and if $h$ is
twice differentiable with bounded derivatives, then
(\ref{SteinGammaEq}) has a solution $f$ which is bounded on
$(-\nu,+\infty)$, differentiable and such that $\|f\|_\infty \leq
\sqrt{2\pi/\nu} \|h\|_\infty$ and $\|f'\|_\infty \leq
\sqrt{2\pi/\nu} \|h'\|_\infty$.

\end{enumerate}
\end{lemme}
Now define
\begin{eqnarray}\label{H1}
\mathscr{H}_1 &=& \{h\in\mathscr{C}^2_b \,:\,\|h\|_\infty
\leq1,\,\|h'\|_\infty\leq 1\}, \label{H22}\\
\mathscr{H}_2 &=& \{h\in\mathscr{C}^2_b\,:\,\|h\|_\infty
\leq1,\,\|h'\|_\infty\leq 1,\, \|h''\|_\infty\leq 1\}, \label{H11}\\
\mathscr{H}_{1,\nu} &=&\mathscr{H}_1\cap\mathscr{C}^2_b(\nu)\\
\label{H2}
\mathscr{H}_{2,\nu} &=& \mathscr{H}_2\cap\mathscr{C}^2_b(\nu)
\end{eqnarray}
where $\mathscr{C}^2_b$ denotes the class of twice differentiable
functions (with support in $\R$) and with bounded derivatives, and
$\mathscr{C}^2_b(\nu)$ denotes the subset of $\mathscr{C}^2_b$
composed of functions with support in $(-\nu,+\infty)$. Note that
point (ii) in the previous statement implies that, by adopting the
notation (\ref{distances}) and for every $\nu>0$ and every real random
variable $Y$ (not necessarily with support in $(-\nu,+\infty)$),
\begin{equation}\label{boundGamma}
d_{\mathscr{H}_{2,\nu}}(Y,F(\nu))\leq
\sup_{f\in\mathscr{F}_{2,\nu}}|E[2(Y+\nu)f'(Y)-Yf(Y)]|
\end{equation}
where $\mathscr{F}_{2,\nu}$ is the class of differentiable
functions with support in $(-\nu,+\infty)$, bounded by 2 and whose first
derivatives are bounded by 1. Analogously, point (iii) implies
that, for every integer $\nu\geq1$,
\begin{equation}\label{boundChi2}
d_{\mathscr{H}_{1,\nu}}(Y,F(\nu))\leq
\sup_{f\in\mathscr{F}_{1,\nu}}|E[2(Y+\nu)f'(Y)-Yf(Y)]|,
\end{equation}
where $\mathscr{F}_{1,\nu}$ is the class of differentiable
functions with support in $(-\nu,+\infty)$, bounded by
$\sqrt{2\pi/\nu}$ and whose first derivatives are also bounded by
$\sqrt{2\pi/\nu}$. A little inspection shows that the following
estimates also hold: for every $\nu>0$ and every random variable
$Y$,
\begin{eqnarray}\label{boundGammabis}
d_{\mathscr{H}_2}(Y,F(\nu))\leq
\sup_{f\in\mathscr{F}_{2}}|E[2(Y+\nu)_+f'(Y)-Yf(Y)]|
\end{eqnarray}
where $\mathscr{F}_{2}$ is the class of functions (defined on
$\R$) that are continuous and differentiable on $\R \backslash
\{\nu\}$, bounded by $\max\{2,2/\nu\}$, and whose first derivatives are
bounded by $\max\{1, 1/\nu+2/\nu^2\}$. Analogously, for every
integer $\nu\geq1$,
\begin{equation}
d_{\mathscr{H}_1}(Y,F(\nu))\leq
\sup_{f\in\mathscr{F}_{1}}|E[2(Y+\nu)_+ f'(Y)-Yf(Y)]|,
\label{boundChi2bis}
\end{equation}
where $\mathscr{F}_{1}$ is the class of functions (on $\R$) that
are continuous and differentiable on $\R \backslash \{\nu\}$,
bounded by $\max\{\sqrt{2\pi/\nu}, 2/\nu\}$, and whose first derivatives
are bounded by $\max\{\sqrt{2\pi/\nu}, 1/\nu+2/\nu^2\}$.

Now, the crucial issue is how to estimate the right-hand side of
(\ref{BoundKol})--(\ref{BoundWass}) and
(\ref{boundGamma})--(\ref{boundChi2bis}) for a given choice of
$Y$. Since Stein's initial contribution \cite{Stein_orig}, an
impressive panoply of techniques has been developed in this
direction (see again \cite{Chen_Shao_AOP} or \cite{Reinert_sur}
for a survey; here, we shall quote e.g.: exchangeable pairs
\cite{Stein_book}, diffusion generators \cite{Barbour1990,Goetze},
size-bias transforms \cite{Gold_Rein_97}, zero-bias transforms
\cite{GoldRin_96}, local dependency graphs \cite{Chen_Shao_AOP}
and graphical-geometric rules \cite{Chatterjee_AOP}). Starting
from the next section, we will show that, when working within the
framework of functionals of Gaussian fields, one can very
effectively estimate expressions such as
(\ref{BoundKol})--(\ref{BoundWass}), (\ref{boundGamma}) and
(\ref{boundChi2}) by using techniques of Malliavin calculus.
Interestingly, a central role is played by an infinite dimensional
version of the same \textsl{integration by parts formula} that is
at the very heart of Stein's characterization of the Gaussian
distribution.
\subsection{The basic approach (with some examples)}
Let $\HH$ be a real separable Hilbert space and,
for $q\geq 1$, let $\HH^{\otimes q}$ (resp. $\HH^{\odot q}$) be
the $q$th tensor product (resp. $q$th symmetric tensor product) of
$\HH$.
We write
\begin{equation} \label{ISO!}
X = \{X(h):h \in \HH\}
\end{equation}
to indicate a centered isonormal Gaussian process on $\HH$. For
every $q\geq 1$, we denote by $I_q$ the isometry between
$\HH^{\odot q}$ (equipped with the norm
$\sqrt{q!}\|\cdot\|_{\HH^{\otimes q}}$) and the $q$th Wiener chaos
of $X$. Note that, if $\HH$ is a $\sigma$-finite measure space
with no atoms, then each random variable $I_q(h)$, $h \in
\HH^{\odot q}$, has the form of a multiple Wiener-It\^{o} integral
of order $q$. We denote by ${\rm L^2}(X)={\rm
L}^2(\Omega,\sigma(X),P)$ the space of square integrable
functionals of $X$, and by $\sk^{1,2}$ the domain of the Malliavin
derivative operator $D$ (see the forthcoming Section \ref{SS :
Elements Mall} for precise definitions). Recall that, for every
$F\in\sk^{1,2}$, the derivative $DF$ is a random element with
values in $\HH$.

We start by observing that, thanks to (\ref{SteinTruc}), for every
$h\in\HH$ such that $\|h\|_\HH=1$ and for every smooth function
$f$, we have $E[X(h)f(X(h))]$ $=$ $E[f'(X(h))]$. Our point is that
this last relation is a very particular case of the following
corollary of the celebrated \textsl{integration by parts formula}
of Malliavin calculus: for every $Y\in\sk^{1,2}$ with zero mean,
\begin{equation}\label{clef}
E[Yf(Y)] = E[\langle DY, -DL^{-1}Y\rangle_\HH f'(Y)],
\end{equation}
where the linear operator $L^{-1}$ is the inverse of the
\textsl{generator of the Ornstein-Uhlebeck semigroup}, noted $L$.
The reader is referred to Section \ref{SS : Elements Mall} and
Section \ref{SS : Stein + Parts} for definitions and for a full
discussion of this point; here, we shall note that $L$ is an
infinite-dimensional version of the generator associated with
Ornstein-Uhlenbeck diffusions (see \cite[Section 1.4]{Nbook} for a
proof of this fact), an object which is also crucial in the
Barbour-Götze ``generator approach'' to Stein's method
\cite{Barbour1990, Goetze}.

It follows that, for every $Y\in\sk^{1,2}$ with zero mean, the expressions
appearing on the right-hand side of (\ref{BoundKol})--(\ref{BoundWass}) (or
(\ref{boundGamma})--(\ref{boundChi2bis})) can be assessed by first
replacing $Yf(Y)$ with $\langle DY, -DL^{-1}Y\rangle_\HH f'(Y)$
inside the expectation, and then by evaluating the $\rm L^2$
distance between $1$ (resp. $2Y+2\nu$) and the inner product
$\langle DY, - DL^{-1}Y\rangle_\HH$. In general, these computations
are carried out by first resorting to the representation of
$\langle DY, -DL^{-1}Y\rangle_\HH$ as a (possibly infinite) series
of multiple stochastic integrals. We will see that, when
$Y=I_q(g)$, for $q\geq2$ and some $g\in\HH^{\odot q}$, then
$\langle DY, -DL^{-1}Y\rangle_\HH=q^{-1}\|DY\|^2_\HH$. In
particular, by using this last relation one can deduce bounds
involving quantities that are intimately related to the central
and non-central limit theorems recently proved in \cite{NP},
\cite{NO} and \cite{NouPec07}.

\begin{rem}
{\rm
\begin{enumerate}
\item The crucial equality $E[I_q(g)f(I_q(g))]$ $=$
$E[q^{-1}\|DI_q(g)\|^2_\HH f'(I_q(g))]$, in the case where $f$ is
a complex exponential, has been first used in \cite{NO}, in order
to give refinements (as well as alternate proofs) of the main CLTs
in \cite{NP} and \cite{PT}. The same relation has been later
applied in \cite{NouPec07}, where a characterization of
non-central limit theorems for multiple integrals is provided.
Note that neither \cite{NouPec07} nor \cite{NO} are concerned with
Stein's method or, more generally, with bounds on distances
between probability measures.
\item We will see that formula (\ref{clef}) contains as a special
case a result recently proved by Chatterjee \cite[Lemma
5.3]{Chatterjee_ptrf}, in the context of limit theorems for linear
statistics of eigenvalues of random matrices. The connection
between the two results can be established by means of the
well-known \textsl{Mehler's formula} (see e.g. \cite[Section 8.5,
Ch. I]{MallBook} or \cite[Section 1.4]{Nbook}), providing a
mixture-type representation of the infinite-dimensional
Ornstein-Uhlenbeck semigroup. See Remarks \ref{Mehler1} and
\ref{Mehler2} below for a precise discussion of this point. See
e.g. \cite{MeyerOU} for a detailed presentation of the
infinite-dimensional Ornstein-Uhlebeck semigroup.
\item We stress that the random variable $\langle DY,
-DL^{-1}Y\rangle_\HH$ appearing in (\ref{clef}) is in general
\textsl{not} measurable with respect to $\sigma(Y)$. For instance,
if $X$ is taken to be the Gaussian space generated by a standard
Brownian motion $\{W_t : t\geq0\}$ and $Y=I_2(h)$ with $h\in
L^{2}_s([0,1]^2)$, then $D_t Y=2\int_0^1 h(u,t)dW_u$, $t\in[0,1]$,
and $$\langle DY, -DL^{-1}Y\rangle_{{\rm
L}^2([0,1])}=2\,I_2(h\otimes_1 h) + 2 \|h\|^2_{{\rm
L}^2([0,1]^2)}$$ which is, in general, not measurable with respect
to $\sigma(Y)$ (the symbol $h\otimes_1 h$ indicates a contraction
kernel, an object that will be defined in Section 2).
\item Note that (\ref{clef}) also implies the relation
\begin{equation}\label{LectVI}
E[Yf(Y)]=E[\tau(Y)f'(Y)],
\end{equation}
where $\tau(Y)=E[\langle DY, -DL^{-1}Y\rangle_\HH|Y]$. Some
general results for the existence of a real-valued function $\tau$
satisfying (\ref{LectVI}) are contained e.g. in \cite{CPU}. Note
that, in general, it is very hard to find an analytic expression
for $\tau(Y)$, especially when $Y$ is a random variable with a
very complex structure, such as e.g. a multiple Wiener-Itô
integral. On the other hand, we will see that, in many cases, the
random variable $\langle DY, -DL^{-1}Y\rangle_\HH$ is remarkably
tractable and explicit. See the forthcoming Section \ref{S :
Unif}, which is based on \cite[Lecture VI]{Stein_book}, for a
general discussions of equations of the type (\ref{LectVI}). See
Remark \ref{zero-bias} below for a connection with Goldstein and
Reinert's \textsl{zero bias transform} \cite{Gold_Rein_97}.
\item The
reader is referred to \cite{PrivReveillac} for applications of
integration by parts techniques to the Stein-type estimation of
drifts of Gaussian processes. See \cite{Hsu} for a Stein
characterization of Brownian motions on manifolds by means of
integration by parts formulae. See \cite{Dec_Savy} for a
connection between Stein's method and algebras of operators on
configuration spaces.
\end{enumerate}
}
\end{rem}

Before proceeding to a formal discussion, and in order to motivate
the reader, we shall provide two examples of the kind of results
that we will obtain in the subsequent sections. The first
statement involves double Wiener-Itô integrals, that is, random
variables living in the second chaos of $X$. The proof is given in
Section \ref{S : Proofs}.
\begin{thm}\label{quidechire}
Let $(Z_n)_{n\geq 1}$ be a sequence belonging to the second Wiener
chaos of $X$. \\
1. Assume that $E(Z_n^2)\rightarrow 1$
and $E(Z_n^4)\rightarrow 3$ as $n\rightarrow\infty$. Then
$Z_n
\,\,\,{\stackrel{{\rm Law}}{\longrightarrow}}\,\,\,
 Z\sim\mathscr{N}(0,1)$ as $n\rightarrow\infty$. Moreover,
we have:
$$
d_{\rm TV} (Z_n,Z)\leq 2
\sqrt{\frac16\big|E(Z_n^4)-3\big|+\frac{3+E(Z_n^2)}{2}\,\big|E(Z_n^2)-1\big|}.
$$
2. Fix $\nu>0$ and assume that $E(Z_n^2)\rightarrow 2\nu$ and
$E(Z_n^4)-12E(Z_n^3)\rightarrow 12\nu^2-48\nu$ as
$n\rightarrow\infty$. Then, as $n\rightarrow\infty$, $Z_n
\,\,\,{\stackrel{{\rm Law}}{\longrightarrow}}\,\,\,
 F(\nu)$, where $F(\nu)$ has a centered Gamma distribution of parameter $\nu$. Moreover, we have:
\begin{eqnarray*}
&&d_{\mathscr{H}_2}(Z_n,F(\nu))\\ &\leq& \max\{1,1/\nu, 2/\nu^2\}\, \sqrt{
\frac16\big|E(Z_n^4)\!-\!12E(Z_n^3)\!-\!12\nu^2\!+\!48\nu\big|\!
+\!\frac{\big|8-6\nu+E(Z_n^2)\big|}{2}\,\big|E(Z_n^2)\!-\!2\nu\big|
},
\end{eqnarray*}
where  $\mathscr{H}_2$
is defined by (\ref{H11}).
\end{thm}
Note that, in the statement of Theorem \ref{quidechire}, there is
no mention of Malliavin operators (however, these operators will
appear in the general statements presented in Section \ref{SS :
Stein + Parts}). For instance, when applied to the case where $X$
is the isonormal process generated by a fractional Brownian
motion, the first point of Theorem \ref{quidechire} can be used to
derive the following bound for the Kolmogorov distance in the
Breuer-Major CLT associated with quadratic transformations:
\begin{thm}\label{BMthm}
Let $B$ be a fractional Brownian motion with Hurst index
$H\in(0,\nicefrac34)$. We set $$
\sigma^2_H=\frac12\sum_{t\in\Z}\big(
|t+1|^{2H}+|t-1|^{2H}-2|t|^{2H}\big)^2<\infty,
$$
and
$$
Z_n=\frac{1}{\sigma_H \sqrt{n}}\sum_{k=0}^{n-1}
\big(n^{2H}(B_{(k+1)/n}-B_{k/n})^2-1\big), \quad n\geq 1.
$$
Then,
as $n\rightarrow\infty$,
$Z_n
\,\,\,{\stackrel{{\rm Law}}{\longrightarrow}}\,\,\,
 Z\sim\mathscr{N}(0,1)$. Moreover, there exists a constant $c_H$ (depending
uniquely on $H$) such that, for any $n\geq 1$:
\begin{equation}\label{ciacca}
d_{\rm Kol}(Z_n,Z)\leq \frac{c_H}{n^{\frac12\wedge (\frac32-2H)}}.
\end{equation}
\end{thm}
Note that both Theorem \ref{quidechire} and \ref{BMthm} will be
significantly generalized in Section \ref{SS : Stein + Parts} and
Section \ref{S : Breuer} (see, in particular, the forthcoming
Theorems \ref{main-tool}, \ref{gam-thm} and \ref{BM-rev}).

\begin{rem}
{\rm
\begin{enumerate}
\item When $H=\nicefrac12$, then $B$ is a standard Brownian motion
(and therefore has independent increments), and we recover from
the previous result the rate $n^{-1/2}$, that could be also
obtained by applying the Berry-Ess\'een Theorem \ref{T : BE}. This
rate is still valid for $H<\nicefrac12$. But, for
$H>\nicefrac12$, the rate we can prove in the Breuer-Major CLT
is $n^{2H-\frac32}$.
\item To the authors knowledge, Theorem \ref{BMthm} and its
generalizations are the first Berry-Ess\'een bounds ever established
for the Breuer-Major CLT.
\item To keep the length of this paper within limits, we do not
derive the explicit expression of some of the constants (such as
the quantity $c_H$ in formula (\ref{ciacca})) composing our
bounds. As will become clear later on, the exact value of these
quantities can be deduced by a careful bookkeeping of the bounding
constants appearing at the different stages of the proofs.
\end{enumerate}
}
\end{rem}

\subsection{Plan}

The rest of the paper is organized as follows. In Section \ref{SS
: Elements Mall} we recall some notions of Malliavin calculus. In
Section \ref{SS : Stein + Parts} we state and discuss our main
bounds in Stein-type estimates for functionals of Gaussian fields.
Section \ref{S : Breuer} contains an application to the
Breuer-Major CLT. Section \ref{S : Gamma} deals with Gamma-type
approximations. Section \ref{S : Unif} provides a unified
discussion of approximations by means of absolutely continuous
distributions. Proofs and further refinements are collected in
Section \ref{S : Proofs}.

\section{Elements of Malliavin calculus} \label{SS : Elements Mall}
The reader is referred to \cite{Janson} or \cite{Nbook} for any
unexplained notion discussed in this section. As in (\ref{ISO!}),
we denote by $X = \{X(h):h \in \HH\}$ an isonormal Gaussian
process over $\HH$. By definition, $X$ is a centered Gaussian
family indexed by the elements of $\HH$ and such that, for every
$h,g\in\HH$,
\begin{equation}\label{isometry}
E\big[X(h)X(g)\big]=\langle h,g\rangle_\HH.
\end{equation}
As before, we use the notation ${\rm L}^2(X)$ $=$ ${\rm
L}^2(\Omega,\sigma(X),P)$. It is well-known (see again \cite[Ch.
1]{Nbook} or \cite{Janson}) that any random variable $F$ belonging
to ${\rm L}^2(X)$ admits the following chaotic expansion:
\begin{equation}\label{ChaosExpansion}
F=\sum_{q=0}^\infty I_q(f_q),
\end{equation}
where $I_0(f_0):=E[F]$, the series converges in ${\rm L}^2$ and
the symmetric kernels $f_q\in\HH^{\odot q}$, $q\geq1$, are
uniquely determined by $F$. As already discussed, in the
particular case where $\HH={\rm L}^2(A,\mathscr{A},\mu)$, where
$(A,\mathscr{A})$ is a measurable space and $\mu$ is a
$\sigma$-finite and non-atomic measure, one has that $\HH^{\odot
q}= {\rm L}^2_s(A^q,\mathscr{A}^{\otimes q}, \mu^{\otimes q})$ is
the space of symmetric and square integrable functions on $A^q$.
Moreover, for every $f\in\HH^{\odot q}$, the random variable
$I_q(f)$ coincides with the multiple Wiener-It\^o integral (of
order $q$) of $f$ with respect to $X$ (see \cite[Ch. 1]{Nbook}).
Observe that a random variable of the type $I_q(f)$, with
$f\in\HH^{\odot q}$, has finite moments of all orders (see e.g.
\cite[Ch. VI]{Janson}). See again \cite[Ch. 1]{Nbook} or
\cite{Rota Wall} for a connection between multiple Wiener-Itô and
Hermite polynomials. For every $q\geq 0$, we write $J_q$ to
indicate the orthogonal projection operator on the $q$th Wiener
chaos associated with $X$, so that, if $F \in {\rm
L}^2(\Omega,\mathscr{F},P)$  is as in (\ref{ChaosExpansion}), then
$J_q F=I_q(f_q)$ for every $q\geq 0$.

Let $\{e_k,\,k\ge 1\}$ be a complete orthonormal system in $\HH$.
Given $f\in \HH^{\odot p}$ and $g\in\HH^{\odot q}$, for every
$r=0,\ldots,p\wedge q$, the $r$th contraction of $f$ and $g$ is
the element of $\HH^{\otimes (p+q-2r)}$ defined as
\begin{equation}\label{defContr}
f\otimes_r g = \sum_{i_1,\ldots,i_r=1}^\infty \langle
f,e_{i_1}\otimes \ldots\otimes e_{i_r}\rangle_{\HH^{\otimes r}}
\otimes \langle g,e_{i_1}\otimes \ldots\otimes
e_{i_r}\rangle_{\HH^{\otimes r}}.
\end{equation}
Note that, in the particular case where $\HH={\rm
L}^2(A,\mathscr{A},\mu)$ (with $\mu$ non-atomic), one has that
$$
f\otimes_r g = \int_{A^r}
f(t_1,\ldots,t_{p-r},s_1,\ldots,s_r)\,g(t_{p-r+1},\ldots,
t_{p+q-2r},s_1,\ldots,s_r)d\mu(s_1)\ldots d\mu(s_r).
$$
Moreover, $f\otimes_0 g=f\otimes g$ equals the tensor product of
$f$ and $g$ while, for $p=q$, $f\otimes_p g=\langle
f,g\rangle_{\HH^{\otimes p}}$. Note that, in general (and except
for trivial cases), the contraction $f\otimes_r g$ is \textsl{not}
a symmetric element of $\HH^{\otimes (p+q-2r)}$. The canonical
symmetrization of $f\otimes_r g$ is written
$f\widetilde{\otimes}_r g$. We also have the useful
multiplication formula: if $f\in \HH^{\odot p}$ and $g\in\HH^{\odot q}$,
then
\begin{equation}\label{multiplication}
I_p(f)I_q(g)=\sum_{r=0}^{p\wedge q} r!\binom{p}{r}\binom{q}r
I_{p+q-2r}(f\widetilde{\otimes}_r g).
\end{equation}

Let $\mathscr{S}$ be the set of all smooth cylindrical random
variables of the form $$F = g\big(X(\phi_1), \ldots,
X(\phi_n)\big)$$ where $n\geq 1$, $g : \R^n \rightarrow \R$ is a
smooth function with compact support and $\phi_i\in\EuFrak H$. The
Malliavin derivative of $F$ with respect to $X$ is the element of
${\rm L}^2(\Omega, \EuFrak H)$ defined as
$$DF\; =\; \sum_{i =1}^n \frac{\partial g}{\partial x_i}\big(X(\phi_1), \ldots, X(\phi_n)\big)
\phi_i.$$ In particular, $DX(h) = h$ for every $h\in \HH$. By
iteration, one can define the $m$th derivative $D^m F$ (which is
an element of ${\rm L}^2(\Omega, \HH^{\otimes m})$) for every
$m\geq 2$. As usual, for $m\geq 1$, $\sk^{m,2}$ denotes the
closure of $\mathscr{S}$ with respect to the norm $\| \cdot
\|_{m,2}$, defined by the relation
$$\| F\|_{m,2}^2 \; = \; E\lcr F^2\rcr +\sum_{i=1}^m
E\big[ \| D^i F\|_{\HH^{\otimes i}}^2\big].$$ Note that, if $F\neq
0$ and $F$ is equal to a finite sum of multiple Wiener-Itô
integrals, then $F\in\sk^{m,2}$ for every $m\geq1$ and the law of
$F$ admits a density with respect to the Lebesgue measure. The
Malliavin derivative $D$ satisfies the following \textsl{chain
rule}: if $\varphi:\R^n\rightarrow\R$ is in $\mathscr{C}^1_b$
(that is, the collection of continuously differentiable
functions with a bounded derivative) and if
$\{F_i\}_{i=1,\ldots,n}$ is a vector of elements of $\sk^{1,2}$,
then $\varphi(F_1,\ldots,F_n)\in\sk^{1,2}$ and
$$
D\varphi(F_1,\ldots,F_n)=\sum_{i=1}^n
\frac{\partial\varphi}{\partial x_i} (F_1,\ldots, F_n)DF_i.
$$
Observe that the previous formula still holds when $\varphi$ is a
Lipschitz function and the law of $(F_1,\ldots, F_n)$ has a
density with respect to the Lebesgue measure on $\R^n$ (see e.g.
Proposition 1.2.3 in \cite{Nbook}).

We denote by $\delta$ the
adjoint of the operator $D$, also called the \textsl{divergence
operator}. A random element $u \in {\rm L}^{2}(\Omega, \HH)$
belongs to the domain of $\delta$, noted ${\rm Dom}\delta$, if,
and only if, it satisfies
$$
| E\langle D F,u\rangle_{\EuFrak H}|\leq c_u\,\|F\|_{{\rm
L}^2}\quad\mbox{for any }F\in{\mathscr{S}},
$$
where $c_u$ is a constant depending uniquely on $u$. If $u \in
{\rm Dom}\delta $, then the random variable $ \delta(u)$ is
defined by the duality relationship (customarily called
``integration by parts formula''):
\begin{equation}\label{ipp}
E (F \delta(u))=  E \langle D F, u \rangle_{\HH},
\end{equation}
which holds for every $F \in \sk^{1,2}$. One sometimes needs the
following property: for every $F\in\sk^{1,2}$ and every $u\in{\rm
Dom}\delta$ such that $Fu$ and $F\delta(u)+\langle
DF,u\rangle_\HH$ are square integrable, one has that $Fu\in{\rm
Dom}\delta$ and
\begin{equation}\label{dfu}
\delta(Fu)=F\delta(u) - \langle DF,u\rangle_\HH.
\end{equation}

The operator $L$, acting on square integrable random variables of
the type (\ref{ChaosExpansion}), is defined through the projection
operators $\{J_q\}_{q \geq 0}$ as $L=\sum_{q=0}^\infty -qJ_q,$ and
is called the \textsl{infinitesimal generator of the
Ornstein-Uhlenbeck semigroup}. It verifies the following crucial
property: a random variable $F$ is an element of $\rm{Dom}$$ L
\,\,(=\sk^{2,2})$ if, and only if, $F\in{\rm Dom}\delta D$ (i.e.
$F\in\sk^{1,2}$ and $DF\in{\rm Dom}\delta$), and in this case:
$\delta DF=-LF.$  Note that a random variable $F$ as in
(\ref{ChaosExpansion}) is in $\sk^{1,2}$ (resp. $\sk^{2,2}$) if,
and only if,
$$
\sum_{q=1}^\infty q\|f_q\|^2_{\HH^{\odot q}}<\infty \quad
(\txt{resp.} \: \sum_{q=1}^\infty q^2\|f_q\|^2_{\HH^{\odot
q}}<\infty),
$$
and also $E\big[\|DF\|^2_\HH\big]=\sum_{q\geq 1}
q\|f_q\|^2_{\HH^{\odot q}}$. If $\HH={\rm L}^2(A,\mathscr{A},\mu)$
(with $\mu$ non-atomic), then the derivative of a random variable
$F$ as in (\ref{ChaosExpansion}) can be identified with the
element of $L^2(A\times \Omega)$ given by
\begin{equation}\label{dtf}
D_a F=\sum_{q=1}^\infty qI_{q-1}\big(f_q(\cdot,a)\big), \quad a
\in A .
\end{equation}
We also define the operator $L^{-1}$, which is the
\textsl{inverse} of $L$, as follows: for every $F\in{\rm L^2}(X)$ with zero mean,
we set $L^{-1}F$ $=$ $\sum_{q\geq 1}-\frac{1}{q} J_q(F)$. Note that
$L^{-1}$ is an operator with values in $\sk^{2,2}$.
The following Lemma contains two statements: the first one
(formula (\ref{OUmwi})) is an immediate consequence of the
definition of $L$ and of the relation $\delta D=-L$, whereas the
second (formula (\ref{MallMoments})) corresponds to Lemma 2.1 in
\cite{NouPec07}.
\begin{lemme} \label{tec}
Fix an integer $q \geq 2$ and set $F=I_q(f)$, with $f\in\HH^{\odot
q}$. Then,
\begin{equation} \label{OUmwi}
\delta D F=qF.
\end{equation}
Moreover, for every integer $s\geq 0$,
\begin{equation}\label{MallMoments}
E\big(F^s\|DF\|_\HH^2\big)=\frac{q}{s+1} E\big(F^{s+2}\big).
\end{equation}
\end{lemme}

\section{Stein's method and integration by parts on Wiener
space} \label{SS : Stein + Parts}
\subsection{Gaussian approximations}\label{SS : GaussAPP}
Our first result provides explicit bounds for the normal
approximation of random variables that are Malliavin derivable.
Although its proof is quite easy to obtain, the following
statement will be central for the rest of the paper.
\begin{thm} \label{main-tool} Let $Z\sim\mathscr{N}(0,1)$, and let $F\in\sk^{1,2}$ be such that
$E(F)=0$. Then, the following bounds are in order:
\begin{eqnarray}
d_{\rm W}(F,Z)&\leq &E[(1-\langle DF,
-DL^{-1}F\rangle_\HH)^2]^{1/2}, \label{MalBoundWas}\\
d_{\rm FM}(F,Z)&\leq &4E[(1-\langle DF,
-DL^{-1}F\rangle_\HH)^2]^{1/2}. \label{MalBoundFM}
\end{eqnarray}
If, in addition, the law of $F$ is absolutely continuous with
respect to the Lebesgue measure, one has that
\begin{eqnarray}
d_{\rm Kol}(F,Z)&\leq &E[(1-\langle DF,
-DL^{-1}F\rangle_\HH)^2]^{1/2}, \label{MalBoundKol}\\
d_{\rm TV}(F,Z)&\leq& 2E[(1-\langle DF,
-DL^{-1}F\rangle_\HH)^2]^{1/2}. \label{MalBoundTV}
\end{eqnarray}
\end{thm}
{\bf Proof.} Start by observing that one can write
$F=LL^{-1}F=-\delta DL^{-1}F$. Now let $f$ be a real
differentiable function. By using the integration by parts formula
and the fact that $Df(F)=f'(F)DF$ (note that, for this formula to
hold when $f$ is only a.e. differentiable, one needs $F$ to have
an absolutely continuous law, see Proposition 1.2.3 in \cite{Nbook}), we deduce
$$E(Ff(F))=E[f'(F)\langle DF,-DL^{-1}F\rangle_\HH].$$
It follows that $E[f'(F)-Ff(F)]=E(f'(F)(1-\langle
DF,-DL^{-1}F\rangle_\HH))$ so that relations
(\ref{MalBoundWas})--(\ref{MalBoundTV}) can be deduced from
(\ref{BoundKol})--(\ref{BoundWass}) and the Cauchy-Schwarz
inequality. \fin

We shall now prove that the bounds appearing in the statement
of Theorem \ref{main-tool} can be \textsl{explicitly} computed,
whenever $F$ belongs to a fixed Wiener chaos.

\begin{prop}\label{Gauss_Wiener} Let $q\geq 2$ be an integer, and
let $F=I_q(f)$, where $f\in\HH^{\odot q}$. Then, $\langle DF,
-DL^{-1}F\rangle_\HH=q^{-1}\|DF\|^2_\HH$, and
\begin{eqnarray}\label{Wienerbound}
&&E[(1-\langle DF,
-DL^{-1}F\rangle_\HH)^2]=E[(1-q^{-1}\|DF\|^2_\HH)^2] \label{WienerBoundUU}\\
&&=(1-q!\left\Vert f\right\Vert_{\HH^{\otimes q}} ^{2})^2
+q^{2}\sum_{r=1}^{q-1}(2q-2r)!(r-1) !^2\binom{q-1}{r-1}^{4}\|
f\widetilde{\otimes }_{r} f \|^2_{\HH^{\otimes 2(q-r)}} \label{WienerBoundI} \\
&&\leq (1-q!\left\Vert f\right\Vert_{\HH^{\otimes q}} ^{2})^2
+q^{2}\sum_{r=1}^{q-1}(2q-2r)!(r-1)!^2\binom{q-1}{r-1}^{4}\|
f\otimes_{r} f \|^2_{\HH^{\otimes 2(q-r)}}.\label{WienerBoundF}
\end{eqnarray}
\end{prop}
{\bf Proof.} The equality $\langle DF,
-DL^{-1}F\rangle_\HH=q^{-1}\|DF\|^2_\HH$ is an immediate
consequence of the relation $L^{-1}I_q(f)=-q^{-1}I_q(f)$. From the
multiplication formulae between multiple stochastic integrals, see
(\ref{multiplication}), one deduces that
\begin{equation} \label{normDERIV}
\Vert D[I_q(f)] \Vert_{\HH} ^{2}=qq!\left\Vert
f\right\Vert_{\HH^{\otimes q}} ^{2}+q^{2}\sum_{r=1}^{q-1}\left(
r-1\right) !\binom{q-1}{r-1}^{2}I_{2\left( q-r\right) }\left(
f\widetilde{\otimes }_{r}f\right)
\end{equation}%
(see also \cite[Lemma 2]{NO}). We therefore obtain
(\ref{WienerBoundI}) by using the orthogonality and isometric
properties of multiple stochastic integrals. The inequality in
(\ref{WienerBoundF}) is just a consequence of the relation $\|
f\widetilde{\otimes }_{r} f \|_{\HH^{\otimes 2(q-r)}} \leq \|
f\otimes_{r} f \|_{\HH^{\otimes 2(q-r)}}$. \fin

The previous result should be compared with the forthcoming
Theorem \ref{david-joe}, where we collect the main findings of
\cite{NO} and \cite{NP}. In particular, the combination of
Proposition \ref{Gauss_Wiener} and Theorem \ref{david-joe} shows
that, for every (normalized) sequence $\{F_n : n\geq 1\}$ living
in a fixed Wiener chaos, the bounds given in
(\ref{MalBoundWas})--(\ref{MalBoundFM}) are ``tight'' with respect
to the convergence in distribution towards
$Z\sim\mathscr{N}(0,1)$, in the sense that these bounds converge
to zero if, and only if, $F_n$ converges in distribution to $Z$.

\begin{thm}[\cite{NO, NP}]\label{david-joe}
Fix $q\geq 2$, and consider a sequence $\{F_n : n\geq 1\}$ such
that $F_n=I_q(f_n)$, $n\geq 1$, where $f_n\in\HH^{\odot q}$.
Assume moreover that $E[F_n^2]=q!\|f_n\|^2_{\HH^{\otimes
q}}\rightarrow 1$. Then, the following four conditions are
equivalent, as  $n\to\infty$:
\begin{enumerate}
\item[\rm(i)] $F_n$ converges
in distribution to $Z\sim \mathscr{N}(0,1)$;
\item[\rm(ii)] $E[F_n^4]\rightarrow 3$;
\item[\rm(iii)] for every $r=1,...,q-1$, $\|f_n\otimes_r f_n\|_{\HH^{\otimes 2(q-r)}}\rightarrow
0$;
\item[\rm(iv)] $\|DF_n\|_{\HH}^{2}\rightarrow q$ in ${\rm L^2}$.
\end{enumerate}
\end{thm}

The implications (i) $\leftrightarrow$ (ii) $\leftrightarrow$
(iii) have been first proved in $\cite{NP}$ by means of stochastic
calculus techniques. The fact that (iv) is equivalent to either
one of conditions (i)--(iii) is proved in \cite{NO}. Note that
Theorem \ref{main-tool} and Proposition \ref{Gauss_Wiener} above
provide an alternate proof of the implications (iii)
$\rightarrow$ (iv) $\rightarrow$ (i). The implication
(ii) $\rightarrow$ (i) can be seen as a drastic simplification of
the ``method of moments and cumulants'', that is a customary tool
in order to prove limit theorems for functionals of Gaussian
fields (see e.g. \cite{BM, ChaSlud, GS, Major, Sur}). In \cite{PT}
one can find a multidimensional version of Theorem
\ref{david-joe}.

\begin{rem}
{\rm Theorem \ref{david-joe} and its generalizations have been
applied to a variety of frameworks, such as: $p$-variations of
stochastic integrals with respect to Gaussian processes
\cite{BNCP, CNW}, quadratic functionals of bivariate Gaussian
processes \cite{DPY}, self-intersection local times of fractional
Brownian motion \cite{HN}, approximation schemes for scalar
fractional differential equations \cite{NN}, high-frequency CLTs
for random fields on homogeneous spaces \cite{MP, MP2, Pecc},
needlets analysis on the sphere \cite{BKMP}, estimation of
self-similarity orders \cite{TudViens}, power variations of
iterated Brownian motions \cite{NouPec07_ibm}. We expect that the
new bounds proved in Theorem \ref{main-tool} and Proposition
\ref{Gauss_Wiener} will lead to further refinements of these
results. See Section \ref{S : Breuer} and Section \ref{S : Gamma}
for applications and examples.}
\end{rem}

As shown in the following statement, the combination of
Proposition \ref{Gauss_Wiener} and Theorem \ref{david-joe} implies
that, on any fixed Wiener chaos, the Kolmogorov, total variation
and Wasserstein distances \textsl{metrize the convergence in
distribution towards Gaussian random variables}. Other topological
characterizations of the set of laws of random variables belonging
to a fixed sum of Wiener chaoses are discussed in \cite[Ch.
VI]{Janson}.

\begin{cor} Let the assumptions and notation of Theorem
\ref{david-joe} prevail. Then, the fact that $F_n$ converges in
distribution to $Z\sim\mathscr{N}(0,1)$ is equivalent to either
one of the following conditions:
\begin{enumerate}
\item[\rm(a)] $d_{\rm Kol}(F_n,Z)\rightarrow 0$;
\item[\rm(b)] $d_{\rm TV}(F_n,Z)\rightarrow 0$;
\item[\rm(c)] $d_{\rm W}(F_n,Z)\rightarrow 0$.
\end{enumerate}
\end{cor}
{\bf Proof.} If $F_n \stackrel{\rm Law}{\rightarrow} Z$ then, by
Theorem \ref{david-joe}, we have necessarily that $\|DF_n\|_\HH^2
\rightarrow q$ in ${\rm L^2}$. The desired conclusion follows
immediately from relations (\ref{MalBoundWas}),
(\ref{MalBoundKol}), (\ref{MalBoundTV}) and (\ref{WienerBoundUU}).
\fin

Note that the previous result is not trivial, since the topologies
induced by $d_{\rm Kol}$, $d_{\rm TV}$ and $d_{\rm W}$ are
stronger than convergence in distribution.

\begin{rem}\label{Mehler1}
{\rm (\textsl{Mehler's formula and Stein's method, I}). In
\cite[Lemma 5.3]{Chatterjee_ptrf}, Chatterjee has proved the
following result (we use a notation which is slightly different
from the original statement). Let $Y=g(V)$, where
$V=(V_1,...,V_n)$ is a vector of centered i.i.d. standard Gaussian
random variables, and $g : \R^n \rightarrow \R$ is a smooth function
such that: (i) $g$ and its derivatives have subexponential growth
at infinity, (ii) $E(g(V))=0$, and (iii) $E(g(V)^2)=1$. Then, for
any Lipschitz function $f$, one has that
\begin{equation}\label{341}
E[Yf(Y)] =
E[S(V)f'(Y)],
\end{equation}
where, for every $v=(v_1,...,v_n)\in\R^n$,
\begin{equation}\label{Chat}
S(v) = \int_0^1
\frac{1}{2\sqrt{t}}E\left[\sum_{i=1}^n
\frac{\partial g}{\partial
v_i}(v)\frac{\partial g}{\partial
v_i}(\sqrt{t}v+\sqrt{1-t}V)\right]dt,
\end{equation}
so that, for instance, for $Z\sim\mathscr{N}(0,1)$ and by using
(\ref{boundTV}), Lemma \ref{Stein_Lemma_Gauss} (iii), (\ref{dollars}) and Cauchy-Schwarz inequality,
\begin{equation}\label{Chat2}
d_{\rm TV}(Y,Z)\leq 2E[(S(V)-1)^2]^{1/2}.
\end{equation}
We shall prove that (\ref{341}) is a very special case of (\ref{clef}).
Observe first that, without loss of
generality, we can assume that $V_i=X(h_i)$, where $X$ is an
isonormal process over some Hilbert space of the type $\HH={\rm
L}^2(A,\mathscr{A},\mu)$ and $\{h_1,...,h_n\}$ is an orthonormal
system in $\HH$.
Since $Y=g(V_1,\ldots,V_n)$, we have $D_a Y =
\sum_{i=1}^n\frac{\partial g}{\partial x_i}(V)h_i(a)$. On the other
hand, since $Y$ is centered and square integrable, it  admits a
chaotic representation of the form $Y=\sum_{q\geq 1}I_q(\psi_q)$.
This implies in particular that $D_a Y =\sum_{q=1}^\infty
qI_{q-1}(\psi_q(a,\cdot))$. Moreover, one has that
$-L^{-1}Y\!=\!\!\sum_{q\geq 1}\!\frac{1}{q} I_q(\psi_q)$, so that
$-D_a L^{-1}Y=\sum_{q\geq 1} I_{q-1}(\psi_q(a,\cdot))$. Now, let
$T_z$, $z\geq 0$, denote the (infinite dimensional)
\textsl{Ornstein-Uhlenbeck semigroup}, whose action on random
variables $F\in {\rm L^2}(X)$ is given by $T_z(F)=\sum_{q\geq0}
\mathrm{e}^{-qz}J_q(F)$. We can write
\begin{eqnarray}
\int_0^1 \frac{1}{2\sqrt{t}} \,T_{\ln(1/\sqrt{t})} (D_a Y)dt
&=&\int_0^\infty e^{-z}T_z(D_a Y)dz
=\sum_{q\geq 1}\frac1q J_{q-1}(D_a Y) \notag\\
&=&\sum_{q\geq 1} I_{q-1}(\psi_q(a,\cdot))=-D_a L^{-1}Y
\label{rrr}.
\end{eqnarray}
Now recall that \textsl{Mehler's formula} (see e.g. \cite[formula
(1.54)]{Nbook}) implies that, for every function $f$ with
subexponential growth,
$$
T_z(f(V))=E\big[f(e^{-z}v+\sqrt{1-e^{-2z}}V)\big]\big|_{v=V},
\quad z\geq 0.
$$
In particular, by applying this last relation to the partial
derivatives $\frac{\partial g}{\partial v_i}$, $i=1,...,n$, we
deduce from (\ref{rrr}) that
$$
\int_0^1 \frac{1}{2\sqrt{t}}\,T_{\ln(1/\sqrt{t})}(D_a Y) dt=
\sum_{i=1}^n h_i(a) \int_0^1 \frac{1}{2\sqrt{t}}\,
E\big[\frac{\partial g}{\partial
v_i}(\sqrt{t}\,v+\sqrt{1-t}\,V)\big]dt\,\,\big|_{v=V}.
$$
Consequently, (\ref{341}) follows, since
\begin{eqnarray*}
\langle DY,-DL^{-1}Y\rangle_\HH &= &\left\langle \sum_{i=1}^n
\frac{\partial g}{\partial v_i}(V)h_i, \sum_{i=1}^n \,\, \int_0^1
\frac{1}{2\sqrt{t}}\, E\big[\frac{\partial g}{\partial
v_i}(\sqrt{t}\,v+\sqrt{1-t}\,V)\big]dt\,\,\big|_{v=V} \,\, h_i
\right\rangle_\HH\\
&=&S(V).
\end{eqnarray*}
\fin
See also Houdré and Pérez-Abreu \cite{HouPA} for related
computations in an infinite-dimensional setting.
}
\end{rem}
The following result concerns finite sums of multiple integrals.
\begin{prop}\label{P : Gauss Sums}
For $s\geq 2$, fix $s$ integers $2\leq q_1<\ldots<q_s$. Consider a
sequence of the form
$$
Z_n=\sum_{i=1}^s I_{q_i}(f_n^i),\quad n\geq 1,
$$
where $f_n^i\in\HH^{\odot q_i}$. Set
$$
\mathscr{I}=\left\{(i,j,r)\in\{1,\ldots,s\}^2\times\N:\,1\leq
r\leq q_i\wedge q_j\,\,\,
\text{and}\,\,\,(r,q_i,q_j)\neq(q_i,q_i,q_i)\right\}.
$$
Then,
\begin{eqnarray*}
&&E[(1-\langle DZ_n, -DL^{-1}Z_n\rangle_\HH)^2]\leq 2\left(
1-\sum_{i=1}^s q_i! \|f_n^i\|_{\HH^{\otimes q_i}}^2
\right)^2\\
&&\hskip1cm +2s^2 \sum_{(i,j,r)\in\mathscr{I}}
q_i^2(r-1)!^2\binom{q_i-1}{r-1}^2
\binom{q_j-1}{r-1}^2(q_i+q_j-2r)!\\
&&\hskip7cm\times \|f_n^i\otimes_{q_i-r}f_n^i\|_{\HH^{\otimes 2r}}
\|f_n^j\otimes_{q_j-r}f_n^j\|_{\HH^{\otimes 2r}}.
\end{eqnarray*}
In particular, if (as $n\rightarrow\infty$) $E[Z_n^2]=\sum_{i=1}^s
q_i!\|f_n^i\|_{\HH^{\otimes q_i}}^2\longrightarrow 1$ and if, for
any $i=1,\ldots,s$ and $r=1,\ldots,q_i-1$, one has that
$\|f_n^i\otimes_r f_n^i\|_{\HH^{\otimes 2(q_i-r)}} \longrightarrow
0$ , then $Z_n \,\,\,{\stackrel{{\rm Law}}{\longrightarrow}}\,\,\,
Z\sim\mathscr{N}(0,1)$ as $n\rightarrow\infty$, and the
inequalities in Theorem \ref{main-tool} allow to associate bounds
with this convergence.
\end{prop}

\begin{rem}
{\rm \begin{enumerate} \item In principle, by using Proposition
\ref{P : Gauss Sums} it is possible to prove bounds for limit
theorems involving the Gaussian approximation of {\it infinite}
sums of multiple integrals, such as for instance the CLT proved in
\cite[Th. 4]{HN}.
\item Note that, to obtain the convergence result stated in Proposition \ref{P : Gauss Sums},
one does not need to suppose that the quantity
$E[I_q(f_i)^2]=q_i!\|f_n^i\|_{\HH^{\otimes q_i}}^2$ is convergent
for every $i$. One should compare this finding with the CLTs
proved in \cite{PT}, as well as the Gaussian approximations
established in \cite{Pecc}.
\end{enumerate}}
\end{rem}

\noindent {\bf Proof of Proposition \ref{P : Gauss Sums}}. Observe
first that, without loss of generality, we can assume that $X$ is
an isonormal process over some Hilbert space of the type $\HH={\rm
L}^2(A,\mathscr{A},\mu)$. For every $a\in A$, it is immediately
checked that
$$
D_aZ_n=\sum_{i=1}^s q_i I_{q_i-1} \big(f_n^i(\cdot,a)\big)
$$
and
$$
-D_a(L^{-1}Z_n)=D_a\left(\sum_{i=1}^s
\frac1{q_i}I_{q_i}(f_n^i)\right)= \sum_{i=1}^s I_{q_i-1}
\big(f_n^i(\cdot,a)\big).
$$
This yields, using in particular the multiplication formula
(\ref{multiplication}):
\begin{eqnarray*}
&&\langle DZ_n,-DL^{-1}Z_n\rangle_{\HH}\\
&=&\sum_{i,j=1}^s q_i \int_A I_{q_i-1} \big(f_n^i(\cdot,a)\big)I_{q_j-1} \big(f_n^j(\cdot,a)\big)\mu(da)\\
&=&\sum_{i,j=1}^s q_i \sum_{r=0}^{q_i\wedge q_j-1}
r!\binom{q_i-1}{r}
\binom{q_j-1}{r}I_{q_i+q_j-2-2r}\left(\int_A f_n^i(\cdot,a)\otimes_r f_n^j(\cdot,a)\mu(da)\right)\\
&=&\sum_{i,j=1}^s q_i \sum_{r=0}^{q_i\wedge q_j-1}
r!\binom{q_i-1}{r}
\binom{q_j-1}{r}I_{q_i+q_j-2-2r}\left(f_n^i\otimes_{r+1} f_n^j\right)\\
&=&\sum_{i,j=1}^s q_i \sum_{r=1}^{q_i\wedge q_j}
(r-1)!\binom{q_i-1}{r-1}
\binom{q_j-1}{r-1}I_{q_i+q_j-2r}\left(f_n^i\otimes_{r} f_n^j\right)\\
&=&\sum_{i=1}^s q_i!\|f_n^i\|^2_{\HH^{\otimes q_i}} +
\sum_{(i,j,r)\in\mathscr{I}} q_i(r-1)!\binom{q_i-1}{r-1}
\binom{q_j-1}{r-1}I_{q_i+q_j-2r}\left(f_n^i\otimes_{r}
f_n^j\right).
\end{eqnarray*}
Thus, by using (among others) inequalities of the type
$(a_1+\ldots+a_v)^2\leq v(a_1^2+\ldots+a_v^2)$, the isometric
properties of multiple integrals as well $\|f\widetilde{\otimes}_r
g\|\leq \|f\otimes_r g\|$, we obtain
\begin{eqnarray*}
&&E\left([\langle DZ_n,-DL^{-1}Z_n\rangle_{\HH}-1]^2\right)\\
&\leq& 2\left(1-\sum_{i=1}^s q_i!\|f_n^i\|^2_{\HH^{\otimes q_i}}\right)^2 \\
&&\hskip0.2cm+2 E \left(\sum_{(i,j,r)\in\mathscr{I}}
q_i(r-1)!\binom{q_i-1}{r-1}
\binom{q_j-1}{r-1}I_{q_i+q_j-2r}\left(f_n^i\otimes_{r} f_n^j\right)\right)^2\\
&\leq& 2\left(1-\sum_{i=1}^s q_i!\|f_n^i\|^2_{\HH^{\otimes q_i}}\right)^2 \\
&&\hskip0.2cm+2s^2 \sum_{(i,j,r)\in\mathscr{I}}
q_i^2(r-1)!^2\binom{q_i-1}{r-1}^2
\binom{q_j-1}{r-1}^2(q_i+q_j-2r)!
\|f_n^i\otimes_{r} f_n^j\|_{\HH^{\otimes q_i+q_j-2r}}^2\\
&\leq& 2\left(1-\sum_{i=1}^s q_i!\|f_n^i\|^2_{\HH^{\otimes q_i}}\right)^2 \\
&&\hskip0.2cm+2s^2 \sum_{(i,j,r)\in\mathscr{I}}
q_i^2(r-1)!^2\binom{q_i-1}{r-1}^2
\binom{q_j-1}{r-1}^2(q_i+q_j-2r)!\\
&&\hskip7cm\times \|f_n^i\otimes_{q_i-r}f_n^i\|_{\HH^{\otimes 2r}}
\|f_n^j\otimes_{q_j-r}f_n^j\|_{\HH^{\otimes 2r}},
\end{eqnarray*}
the last inequality being a consequence of the (easily verified)
relation
$$
\|f_n^i\otimes_{r} f_n^j\|_{\HH^{\otimes q_i+q_j-2r}} ^2=\langle
f_n^i\otimes_{q_i-r} f_n^i,
f_n^j\otimes_{q_j-r}f_n^j\rangle_{\HH^{\otimes 2r}}.
$$ \fin
\subsection{A property of $\langle DF, -DL^{-1}F\rangle_\HH$}
Before dealing with Gamma approximations, we shall prove the a.s.
positivity of a specific projection of the random variable
$\langle DF, -DL^{-1}F\rangle_\HH$ appearing in Theorem
\ref{main-tool}. This fact will be used in the proof of the main
result of the next section.
\begin{prop}\label{P : positivity}
Let $F\in\sk^{1,2}$. Then, $P$-a.s.,
\begin{equation}\label{POS}
E[\langle DF, -DL^{-1}F\rangle_\HH | F]\geq 0.
\end{equation}
\end{prop}
{\bf Proof.} Let $g$ be a non-negative real function, and set
$G(x)=\int_0^x g(t)dt$, with the usual convention
$\int_0^x=-\int_x^0$ for $x<0$. Since $G$ is increasing and
vanishing at zero, we have $xG(x)\geq 0$ for all $x\in\R$. In
particular, $E(FG(F))\geq 0$. Moreover,
$$ E[F\,G(F)] = E[\langle DF, -DL^{-1}F\rangle_\HH\, g(F)]
= E[E[\langle DF, -DL^{-1}F\rangle_\HH | F]g(F)].$$ We therefore
deduce that
$$E[E[\langle DF, -DL^{-1}F\rangle_\HH |
F]\1_A]\geq 0$$ for any $\sigma(F)$-measurable set $A$. This
implies the desired conclusion. \fin
\begin{rem}\label{zero-bias}
{\rm According to Goldstein and Reinert \cite{Gold_Rein_97}, for
$F$ as in the previous statement, there exists a random variable
$F^*$ having the $F$-\textsl{zero biased} distribution, that is,
$F^*$ is such that, for every absolutely continuous function $f$,
$$ E[f'(F^*)] = E[Ff(F)].$$
By the computations made in the previous proof, one also has that
$$ E[g(F^*)] = E[\langle DF, -DL^{-1}F\rangle_\HH g(F)],$$
for any real-valued and smooth function $g$. This implies, in
particular, that the conditional expectation $E[\langle DF,
-DL^{-1}F\rangle_\HH | F] $ is a version of the Radon-Nikodym
derivative of the law of $F^*$ with respect to the law of $F$,
whenever the two laws are equivalent. }
\end{rem}

\subsection{Gamma approximations}\label{SS : GammaAPP}
We now combine Malliavin calculus with the Gamma approximations
discussed in the second part of Section \ref{SS : Stein Intro}.
\begin{thm}\label{gam-thm}
Fix $\nu>0$ and let $F(\nu)$ have a centered Gamma distribution
with parameter $\nu$. Let $G\in\sk^{1,2}$ be such that $E(G)=0$
and the law of $G$ is absolutely continuous with respect to the
Lebesgue measure. Then:
\begin{equation}
d_{\mathscr{H}_2}(G,F(\nu))\leq K_2 E[(2\nu+2G-\langle DG,
-DL^{-1}G\rangle_\HH)^2]^{1/2}, \label{K1}
\end{equation}
and, if $\nu\geq 1$ is an integer,
\begin{equation}
d_{\mathscr{H}_1}(G,F(\nu))\leq K_1 E[(2\nu+2G-\langle DG,
-DL^{-1}G\rangle_\HH)^2]^{1/2}, \label{K2}
\end{equation}
where $\mathscr{H}_1$ and $\mathscr{H}_2$ are defined in
(\ref{H22})--(\ref{H11}), $K_1$ $ \!\!:= $ $ \!\!\max\{\!\sqrt{2\pi/\nu},\!$
$1/\nu+2/\nu^2\!\}$ and
$K_2$ $:=$ $ \max\{1,1/\nu+ 2/\nu^2\}$.
\end{thm}
{\bf Proof.} We will only prove (\ref{K1}), the proof of
(\ref{K2}) being analogous. Fix $\nu>0$. Thanks to
(\ref{boundGammabis}) and (\ref{clef}) (in the case $Y=G$) and by
applying Cauchy-Schwarz, we deduce that
\begin{eqnarray*}
d_{\mathscr{H}_2}(G,F(\nu))&\leq&
\sup_{\mathscr{F}_2}|E[f'(G)(2(\nu+G)_+ -\langle DG,
-DL^{-1}G\rangle_\HH)]|\\
&\leq & K_2 \times E[(2(\nu+G)_+ - \langle DG,
-DL^{-1}G\rangle_\HH)^2]^{1/2}\\
&\leq & K_2 \times E[(2(\nu+G) - \langle DG,
-DL^{-1}G\rangle_\HH)^2]^{1/2},
\end{eqnarray*}
where the last inequality is a consequence of the fact that
$E[\langle DG, -DL^{-1}G\rangle_\HH|G]$ $\geq0$ (thanks to
Proposition \ref{P : positivity}). \fin

\begin{rem}\label{Mehler2}
{\rm
(\textsl{Mehler's formula and Stein's method, II}). Define
$Y=g(V)$ as in Remark \ref{Mehler1}. Then, since (\ref{Chat}) and
(\ref{Chat2}) are in order, one deduces from Theorem \ref{gam-thm}
that, for every $\nu>0$, $$ d_{\mathscr{H}_2} (Y,F(\nu))\leq K_2
E[(2\nu+2Y-S(V))^2]^{1/2}.$$ An analogous estimate holds for
$d_{\mathscr{H}_1}$, when applied to the case where $\nu\geq1$ is
an integer.
}
\end{rem}
We will now connect the previous results to the main findings of
\cite{NouPec07}. To do this, we shall provide explicit estimates
of the bounds appearing in Theorem \ref{gam-thm}, in the case
where $G$ belongs to a fixed Wiener chaos of \textsl{even order}
$q$.

\begin{prop}\label{Gamma_Wiener} Let $q\geq 2$ be an even integer, and
let $G=I_q(g)$, where $g\in\HH^{\odot q}$. Then,
\begin{eqnarray}\label{WienerboundG}
&&E[(2\nu+2G - \langle DG,
-DL^{-1}G\rangle_\HH)^2] =E[(2\nu+2G-q^{-1}\|DG\|^2_\HH)^2]\\
&& \leq (2\nu-q!\left\Vert g\right\Vert_{\HH^{\otimes q}} ^{2})^2+ \nonumber\\
&& \quad \quad \quad
+q^{2}\sum_{\stackrel{r\in\{1,...,q-1\}}{r\neq q/2}}(2q-2r)!(r-1)
!^2\binom{q-1}{r-1}^{4}\| g\otimes_{r} g
\|^2_{\HH^{\otimes 2(q-r)}} + \nonumber \\
&&\quad\quad\quad\quad\quad\quad\quad\quad\quad\quad\quad\quad\quad\quad\quad\quad\quad\quad\quad\quad+
4q!\left\|c_q^{-1}\times g\widetilde{\otimes}_{q/2} g-g
\right\|^2_{\HH^{\otimes q}}, \nonumber
\end{eqnarray}
where
\begin{equation}\label{ciqqu}
c_{q}:=\frac{1}{\left(q/2\right)!\binom{q-1}{q/2-1}^{2}}=\frac{4}{\left(
q/2\right)!\binom{q}{q/2}^{2}}.
\end{equation}
\end{prop}
{\bf Proof.} By using (\ref{normDERIV}) we deduce that
\begin{eqnarray*}
q^{-1}\Vert DG \Vert_{\HH} ^{2}-2\nu-2G&=&(q!\left\Vert
g\right\Vert_{\HH^{\otimes q}} ^{2}-2\nu) +\\
&+&q\sum_{\stackrel{r\in\{1,...,q-1\}}{r\neq q/2}}\left(
r-1\right) !\binom{q-1}{r-1}^{2}I_{2\left( q-r\right) }\left(
g\widetilde{\otimes }_{r}g\right)+\\
&&\quad \quad \quad \quad \quad
+q(q/2-1)!\binom{q-1}{q/2-1}I_q(g\widetilde{\otimes }_{q/2}g-2g).
\end{eqnarray*}
The conclusion is obtained by using the isometric properties of
multiple Wiener-Itô integrals, as well as the relation $\|
g\widetilde{\otimes }_{r} g \|_{\HH^{\otimes 2(q-r)}} \leq \|
g\otimes_{r} g \|_{\HH^{\otimes 2(q-r)}}$, for every
$r\in\{1,...,q-1\}$ such that $r\neq q/2$. \fin By using
Proposition \ref{Gamma_Wiener}, we immediately recover the
implications (iv) $\rightarrow$ (iii) $\rightarrow$ (i) in the
statement of the following result, recently proved in \cite[Th.
1.2]{NouPec07}.
\begin{thm}[\cite{NouPec07}]\label{ivan-joe}
Let $\nu>0$ and let $F(\nu)$ have a centered Gamma distribution
with parameter $\nu$. Fix an {\rm even} integer $q\geq 2$, and
define $c_{q}$ according to (\ref{ciqqu}). Consider a sequence of
the type $G_n=I_q(g_n)$, where $n\geq1$ and $g_n\in\HH^{\odot q}$,
and suppose that
$$
\lim_{n\rightarrow\infty}
E\big[G_n^2\big]=\lim_{n\rightarrow\infty}
q!\|g_n\|^2_{\HH^{\otimes q}}=2\nu.
$$
Then, the following four conditions are equivalent:
\begin{enumerate}
\item[\rm(i)] as $n\to\infty$, the sequence $(G_n)_{n\geq 1}$ converges  in distribution to
$F(\nu)$;
\item[\rm(ii)] $\lim_{n\to\infty}E[G_n^4]-12E[G_n^3]=12\nu^2-48\nu$;
\item[\rm(iii)] as $n\to\infty$, $\|DG_n\|_{\HH}^{2} - 2q G_n \longrightarrow 2q\nu$ in ${\rm L}^2$.
\item[\rm(iv)] $\lim_{n\to\infty}\|g_n\widetilde{\otimes}_{q/2} g_n-c_q\times g_n\|_{\HH^{\otimes q}}=0,$
\, where $c_q$ is given by (\ref{ciqqu}), and \,
$\lim_{n\to\infty}\|g_n\otimes_{r} g_n\|_{\HH^{\otimes
2(q-r)}}=0$, for every $r=1,...,q-1$ such that $r\neq q/2$.
\end{enumerate}
\end{thm}

\smallskip

\noindent Observe that $E(F(\nu)^2)=2\nu$, $E(F(\nu)^3)=8\nu$ and
$E(F(\nu)^4)=48\nu+12\nu^2$, so that the implication (ii)
$\rightarrow$ (i) in the previous statement can be seen as a
further simplification of the \textsl{method of moments and
cumulants}, as applied to non-central limit theorems (see e.g.
\cite{SUr2}, and the references therein, for a survey of classic
non-central limit theorems). Also, the combination of Proposition
\ref{Gamma_Wiener} and Theorem \ref{ivan-joe} shows that, inside a
fixed Wiener chaos of even order, one has that: (i)
$d_{\mathscr{H}_2}$ metrizes the weak convergence towards centered
Gamma distributions, and (ii) $d_{\mathscr{H}_1}$ metrizes the
weak convergence towards centered $\chi^2$ distributions with
arbitrary degrees of freedom.

The following result concerns the Gamma approximation of a sum of
two multiple integrals. Note, at the cost of a quite heavy
notation, one could easily establish analogous estimates for sums
of three or more integrals. The reader should compare this result
with Proposition \ref{P : Gauss Sums}.

\begin{prop}\label{P : Gamma Sums}
Fix two real numbers $\nu_1,\nu_2 > 0$, as well as two even
integers $2\leq q_1< q_2$. Set $\nu=\nu_1+\nu_2$ and suppose (for
the sake of simplicity) that $q_2> 2q_1$. Consider a sequence of
the form
$$
Z_n=I_{q_1}(f_n^1)+I_{q_2}(f_n^2),\quad n\geq 1,
$$
where $f_n^i\in\HH^{\odot q_i}$. Set
$$
\mathscr{J}=\left\{(i,j,r)\in\{1,2\}^2\times\N:\,1\leq r\leq
q_i\wedge q_j\,\mbox{and, whenever $i=j$, $r\neq q_i\mbox{ and
}r\neq \frac{q_i}2$} \right\}.
$$
Then
\begin{eqnarray}
&&E[(2Z_n+2\nu-\langle DZ_n,
-DL^{-1}Z_n\rangle_\HH)^2] \notag\\
&\leq& 3\left(2\nu-\sum_{i=1,2} q_i!\|f_n^i\|^2_{\HH^{\otimes
q_i}}\right)^2 +24\sum_{i=1,2}c_{q_i}^{-2}q_i! \,\|
f_n^i\widetilde\otimes_{\nicefrac{q_i}2}f_n^i
-c_{q_i}\times f_n^i\|_{\HH^{\otimes q_i}}^2 \notag \\
&&\hskip0.2cm+12 \sum_{(i,j,r)\in\mathscr{J}}
q_i^2(r-1)!^2\binom{q_i-1}{r-1}^2
\binom{q_j-1}{r-1}^2(q_i+q_j-2r)! \label{cool} \\
&&\hskip7cm\times \|f_n^i\otimes_{q_i-r}f_n^i\|_{\HH^{\otimes 2r}}
\|f_n^j\otimes_{q_j-r}f_n^j\|_{\HH^{\otimes 2r}}. \notag
\end{eqnarray}
In particular, if
\begin{enumerate}
\item[(i)] $E[Z_n^2]=\sum_{i=1,2} q_i!\|f_n^i\|_{\HH^{\otimes q_i}}^2\longrightarrow 2\nu$ as $n\rightarrow\infty$,
\item[(ii)] for $i=1,2$,
$\|f_n^i\widetilde\otimes_{\nicefrac{q_i}2} f_n^i - c_{q_i}\times
f_n^i\|_{\HH^{\otimes q_i}}\longrightarrow 0$ as
$n\rightarrow\infty$, where $c_{q_i}$ is defined in Theorem
\ref{ivan-joe},
\item[(iii)] for any $i=1,2$ and $r=1,\ldots,q_i-1$ such that $r\neq \frac{q_i}2$,
$\|f_n^i\otimes_r f_n^i\|_{\HH^{\otimes 2(q_i-r)}}\longrightarrow
0$ as $n\rightarrow\infty$,
\end{enumerate}
then $Z_n \,\,\,{\stackrel{{\rm Law}}{\longrightarrow}}\,\,\,
F(\nu)$ as $n\rightarrow\infty$, and the combination of Theorem
\ref{main-tool} and (\ref{cool}) allows to associate explicit
bounds with this convergence.
\end{prop}
{\bf Proof of Proposition \ref{P : Gamma Sums}}. We have (see the
proof of Proposition \ref{P : Gauss Sums})
\begin{eqnarray*}
&&\langle DZ_n,-DL^{-1}Z_n\rangle_{\HH}-2Z_n-2\nu\\
&&=\left(\sum_{i=1,2}q_i!\|f_n^i\|^2_{\HH^{\otimes
q_i}}-2\nu\right)
+\sum_{i=1,2}2\,c_{q_i}^{-1}\,I_{q_i}(f_n^i\widetilde\otimes_{\nicefrac{q_i}2}f_n^i
-c_{q_i}\times f_n^i)\\
&&\hskip1cm+\sum_{(i,j,r)\in\mathscr{J}}
q_i(r-1)!\binom{q_i-1}{r-1}
\binom{q_j-1}{r-1}I_{q_i+q_j-2r}\left(f_n^i\otimes_{r}
f_n^j\right).
\end{eqnarray*}
Thus
\begin{eqnarray*}
&&E\left([\langle DZ_n,-DL^{-1}Z_n\rangle_{\HH}-2Z_n-2\nu]^2\right)\\
&\leq& 3\left(2\nu-\sum_{i=1,2} q_i!\|f_n^i\|^2_{\HH^{\otimes
q_i}}\right)^2 +24\sum_{i=1,2}c_{q_i}^{-2}q_i! \,\|
f_n^i\widetilde\otimes_{\nicefrac{q_i}2}f_n^i
-c_{q_i}\times f_n^i\|_{\HH^{\otimes q_i}}^2\\
&&\hskip0.2cm+3\, E \left(\sum_{(i,j,r)\in\mathscr{J}}
q_i(r-1)!\binom{q_i-1}{r-1}
\binom{q_j-1}{r-1}I_{q_i+q_j-2r}\left(f_n^i\otimes_{r} f_n^j\right)\right)^2\\
&\leq& 3\left(2\nu-\sum_{i=1,2} q_i!\|f_n^i\|^2_{\HH^{\otimes
q_i}}\right)^2 +24\sum_{i=1,2}c_{q_i}^{-2}q_i! \,\|
f_n^i\widetilde\otimes_{\nicefrac{q_i}2}f_n^i
-c_{q_i}\times f_n^i\|_{\HH^{\otimes q_i}}^2\\
&&\hskip0.2cm+12 \sum_{(i,j,r)\in\mathscr{J}}
q_i^2(r-1)!^2\binom{q_i-1}{r-1}^2
\binom{q_j-1}{r-1}^2(q_i+q_j-2r)!
\|f_n^i\otimes_{r} f_n^j\|_{\HH^{\otimes q_i+q_j-2r}}^2\\
&\leq& 3\left(2\nu-\sum_{i=1,2} q_i!\|f_n^i\|^2_{\HH^{\otimes
q_i}}\right)^2 +24\sum_{i=1,2}c_{q_i}^{-2}q_i! \,\|
f_n^i\widetilde\otimes_{\nicefrac{q_i}2}f_n^i
-c_{q_i}\times f_n^i\|_{\HH^{\otimes q_i}}^2\\
&&\hskip0.2cm+12 \sum_{(i,j,r)\in\mathscr{J}}
q_i^2(r-1)!^2\binom{q_i-1}{r-1}^2
\binom{q_j-1}{r-1}^2(q_i+q_j-2r)!\\
&&\hskip7cm\times \|f_n^i\otimes_{q_i-r}f_n^i\|_{\HH^{\otimes 2r}}
\|f_n^j\otimes_{q_j-r}f_n^j\|_{\HH^{\otimes 2r}}.
\end{eqnarray*} \fin

\section{Berry-Ess\'een bounds in the Breuer-Major CLT} \label{S :
Breuer}

In this section, we use our main results in order to derive an
explicit Berry-Ess\'een bound for the celebrated
\textsl{Breuer-Major CLT} for Gaussian-subordinated random
sequences. For simplicity, we focus on sequences that can be
represented as Hermite-type functions of the (normalized)
increments of a fractional Brownian motion. Our framework include
examples of Gaussian sequences whose autocovariance functions
display long dependence. Plainly, the techniques developed in this
paper can also accommodate the analysis of more general
transformations (for instance, obtained from functions with an
arbitrary \textsl{Hermite rank} -- see \cite{T1}), as well as
alternative covariance structures.
\subsection{General setup}
We recall that a \textsl{fractional Brownian motion} (fBm)
$B=\{B_t:t\in[0,1]\}$, of Hurst
index $H\in(0,1)$, is a centered Gaussian process, started from zero and with covariance function $E%
(B_{s}B_{t})=R_{H}(s,t)$, where
$$
R_{H}(s,t)=\frac{1}{2}\left( t^{2H}+s^{2H}-|t-s|^{2H}\right);\quad
s,t\in[0,1].
$$
If $H=\nicefrac12$, then $R_{H}(s,t)=\min(s,t)$ and $B$ is a
standard Brownian motion. For any choice of the Hurst parameter
$H\in(0,1)$, the Gaussian space generated by $B$ can be identified
with an isonormal Gaussian process of the type
$X=\{X(h):h\in\HH\}$, where the real and separable Hilbert space
$\EuFrak H$ is defined as follows: (i) denote by $\mathscr{E}$ the
set of all $\mathbb{R}$-valued step functions on $[0,1]$, (ii)
define $\EuFrak H$ as the Hilbert space obtained by closing
$\mathscr{E}$ with respect to the scalar product
$$
\left\langle
{\mathbf{1}}_{[0,t]},{\mathbf{1}}_{[0,s]}\right\rangle _{\EuFrak
H}=R_{H}(t,s).
$$
In particular, with such a notation one has that
$B_t=X(\mathbf{1}_{[0,t]})$. Note that, if $H=\nicefrac12$, then
$\HH=L^{2}[0,1]$; when $H>\nicefrac12$, the space $\HH$ coincides
with the space of distributions $f$ such that
$s^{\frac12-H}\mathscr{I}_{0+}^{H-\frac12}(f(u)u^{H-\frac12})(s)$
belongs to $L^{2}[0,1]$; when $H<\nicefrac12$ one has that $\HH$
is $\mathscr{I}_{0+}^{H-\frac{1}{2}}(L^{2}[0,1])$. Here,
$\mathscr{I}_{0+}^{H-\frac{1}{2}}$ denotes the action of the
\textsl{fractional Riemann-Liouville operator}, defined as
$$
\mathscr{I}_{0+}^{H-\frac12} f(x)=\frac1{\Gamma(H-\frac
12)}\int_0^x (x-y)^{H-\frac32}f(y)dy.
$$
The reader is referred e.g. to \cite%
{Nbook} for more details on fBm and fractional operators.
\subsection{A Berry-Ess\'een bound}
In what follows, we will be interested in the
asymptotic behaviour (as $n\rightarrow \infty$) of random vectors
that are subordinated to the array
\begin{equation}\label{fBmInc}
V_{n,H}=\{n^H(B_{(k+1)/n}-B_{k/n}): k=0,...,n-1\},\:\: n\geq1.
\end{equation}
Note that, for every $n\geq 1$, the law of $V_{n,H}$ in
(\ref{fBmInc}) coincides with the law of the first $n$ instants of
a centered stationary Gaussian sequence indexed by $\{0,1,2,...\}$
and with autocovariance function given by
$$\rho_H(k)=\frac12(|k+1|^{2H}+|k-1|^{2H}-2|k|^{2H}),\quad
k\in\Z$$ (in particular, $\rho_H(0)=1$ and
$\rho_H(k)=\rho_H(-k)$). From this last expression, one deduces
that the components of the vector $V_{n,H}$ are: (a) i.i.d. for
$H=\nicefrac12$, (b) negatively correlated for
$H\in(0,\nicefrac12)$ and (c) positively correlated for
$H\in(\nicefrac12,1)$. In particular, if $H\in(\nicefrac12,1)$,
then $\sum_k \rho_H(k)=+\infty$: in this case, one customarily
says that $\rho_H$ exhibits \textsl{long-range dependence} (or,
equivalently, \textsl{long memory} -- see e.g. \cite{T3} for a
general discussion of this point).

\smallskip

\noindent Now denote by $H_q$, $q\geq 2$, the $q$th Hermite
polynomial, defined as
$$H_q(x)=\frac{(-1)^q}{q!}e^{\frac{x^2}{2}} \frac{d^q}{dx^q}e^{-\frac{x^2}{2}}, \; x\in\R.$$
For instance, $H_2(x)=(x^2-1)/2$, $H_3(x)=(x^3-3x)/6$, and so on.
Finally, set
$$\sigma=\sqrt{\frac1{q!}\sum_{t\in\Z}\rho_H(t)^q},$$
and define
\begin{equation}\label{Z_n}
Z_n=\frac{1}{\sigma\sqrt{n}}\sum_{k=0}^{n-1}
H_q\big(n^H(B_{(k+1)/n}-B_{k/n})\big)
=\frac{n^{qH-\frac12}}{q!\sigma}\sum_{k=0}^{n-1}
I_q(\delta_{k/n}^{\otimes q}),
\end{equation}
where $I_q$ denotes the $q$th multiple integral with respect to
the isonormal process associated with $B$ (see Section \ref{SS :
Elements Mall}). For simplicity, here (and for the rest of this
section) we write $\delta_{k/n}$ instead of ${\bf
1}_{[k/n,(k+1)/n]}$, and also $\delta_{k/n}^{\otimes
q}=\delta_{k/n}\otimes\cdot\cdot\cdot \otimes\delta_{k/n}$ ($q$ times). Note
that in (\ref{Z_n}) we have used the standard relation: $q!H_q(h)
= I_q(h^{\otimes q})$ for every $h\in\HH$ such that $\|h\|_\HH=1$
(see e.g. \cite[Ch. 1]{Nbook}).

Now observe that, for every $q\geq2$, one has that
$\sum_t|\rho_H(t)|^q<\infty$
if, and only if, $H\in(0,\frac{2q-1}{2q})$.
Moreover, in this case, $E(Z^2_n)\rightarrow 1$ as $n\rightarrow\infty$. As a consequence,
according to Breuer and Major's well-known result \cite[Theorem
1]{BM}, as $n\rightarrow\infty$
$$
Z_n\rightarrow Z\sim\mathscr{N}(0,1)\quad \text{in distribution.}
$$
To the authors' knowledge, the following statement contains the
first Berry-Ess\'een bound ever proved for the Breuer-Major CLT:
\begin{thm}\label{BM-rev}
As $n\rightarrow\infty$, $Z_n$ converges in law towards
$Z\sim\mathscr{N}(0,1)$. Moreover, there exists a constant $c_H$,
depending uniquely on $H$, such that, for any $n\geq 1$:
\begin{eqnarray*}
\sup_{z\in\R}|P(Z_n\leq z)-P(Z\leq z)| &\leq& c_H \times
\left\{\begin{array}{lll}
n^{-\frac12}&\quad\mbox{if $H\in (0,\frac12]$}\\
\\
n^{H-1}&\quad\mbox{if $H\in [\frac12,\frac{2q-3}{2q-2}]$}\\
\\
n^{qH-q+\frac12}&\quad\mbox{if $H\in [\frac{2q-3}{2q-2},\frac{2q-1}{2q})$}\\
\end{array}\right.
\end{eqnarray*}
\end{thm}
\begin{rem}
{\rm
\begin{enumerate}
\item Theorem \ref{BMthm} (see the Introduction) can be proved by simply setting $q=2$ in Theorem \ref{BM-rev}.
Observe that in this case one has $\frac{2q-3}{2q-2}=\frac12$, so
that the middle line in the previous display becomes immaterial.
\item When $H>\frac{2q-1}{2q}$, the sequence $Z_n$ does not
converge in law towards a Gaussian random variable. Indeed, in
this case a non-central limit theorem takes place. See Breton and
Nourdin \cite{BreNou} for bounds associated with this convergence.
\item As discussed in \cite[p. 429]{BM}, it is in general not possible to derive CLTs such as the one
in Theorem \ref{BM-rev} from mixing-type conditions. In
particular, it seems unfeasible to deduce Theorem \ref{BM-rev}
from any mixing characterization of the increments of fractional
Brownian motion (as the one proved e.g. by Picard in \cite[Theorem
A.1]{picard}). See e.g. Tikhomirov \cite{Tikh} for general
derivations of Berry-Ess\'een bounds from strong mixing conditions.
\end{enumerate}
}
\end{rem}
\subsection{Proof of Theorem \ref{BM-rev}}
We have
$$
DZ_n=\frac{n^{qH-\frac12}}{(q-1)!\sigma}\sum_{k=0}^{n-1}
I_{q-1}(\delta_{k/n}^{\otimes q-1})\delta_{k/n},
$$
hence
$$
\|DZ_n\|^2_\HH=\frac{n^{2qH-1}}{(q-1)!^2\sigma^2}\sum_{k,\ell=0}^{n-1}
I_{q-1}(\delta_{k/n}^{\otimes q-1})
I_{q-1}(\delta_{\ell/n}^{\otimes q-1})
\langle \delta_{k/n},\delta_{\ell/n}\rangle_\HH.
$$
By the multiplication formula (\ref{multiplication}):
$$
I_{q-1}(\delta_{k/n}^{\otimes q-1})
I_{q-1}(\delta_{\ell/n}^{\otimes q-1})
=\sum_{r=0}^{q-1}
r!\binom{q-1}{r}^2\,
I_{2q-2-2r}\big(
\delta_{k/n}^{\otimes q-1-r}\widetilde{\otimes}\delta_{\ell/n}^{q-1-r}
\big)
\langle\delta_{k/n},\delta_{\ell/n}\rangle_\HH^r.
$$
Consequently,
$$
\|DZ_n\|^2_\HH=\frac{n^{2qH-1}}{(q-1)!^2\sigma^2}
\sum_{r=0}^{q-1}r!\binom{q-1}{r}^2\,
\sum_{k,\ell=0}^{n-1}
I_{2q-2-2r}\big(
\delta_{k/n}^{\otimes q-1-r}\widetilde{\otimes}\delta_{\ell/n}^{q-1-r}
\big)
\langle\delta_{k/n},\delta_{\ell/n}\rangle_\HH^{r+1}.
$$
Thus, we can write
\begin{eqnarray*}
\frac1q\|DZ_n\|^2_\HH -1 = \sum_{r=0}^{q-1} A_r(n) - 1
\end{eqnarray*}
where
\begin{eqnarray*}
A_r(n)=\frac{r!\binom{q-1}{r}^2}{q(q-1)!^2\sigma^2}n^{2qH-1}
\sum_{k,\ell=0}^{n-1}
I_{2q-2-2r}\big(
\delta_{k/n}^{\otimes q-1-r}\widetilde{\otimes}\delta_{\ell/n}^{q-1-r}
\big)
\langle\delta_{k/n},\delta_{\ell/n}\rangle_\HH^{r+1}.
\end{eqnarray*}
We will need the following easy Lemma (the proof is omitted). Here
and for the rest of the proof of Theorem \ref{BM-rev}, the
notation $a_n\trianglelefteqslant b_n$ means that $\sup_{n\geq 1}
|a_n|/|b_n|<\infty$.
\begin{lemme}\label{lm-int}
1. We have $\rho_H(n)\trianglelefteqslant |n|^{2H-2}$.\\
2. For any $\alpha\in\R$, we have
$$
\sum_{k=1}^{n-1} k^\alpha \trianglelefteqslant 1+n^{\alpha+1}.
$$
3. If $\alpha\in(-\infty,-1)$, we have
$$
\sum_{k= n}^{\infty} k^\alpha \trianglelefteqslant n^{\alpha+1}.
$$
\end{lemme}
By using elementary computations (in particular, observe that
$n^{2H}\langle\delta_{k/n},\delta_{\ell/n}\rangle_\HH =
\rho_H(k-\ell)$) and then Lemma \ref{lm-int}, it is easy to check
that
\begin{eqnarray*}
A_{q-1}(n)-1&=&\frac{1}{q!\sigma^2}n^{2qH-1}\sum_{k,\ell=0}^{n-1}\langle
\delta_{k/n},\delta_{\ell/n}\rangle_\HH^q-1\\
&=&\frac{1}{q!\sigma^2}\left(\frac1n\sum_{k,\ell=0}^{n-1}
\rho_H(k-\ell)^q - \sum_{t\in\Z}\rho_H(t)^q\right)\\
&=&\frac{1}{q!\sigma^2}\left(
\frac1n\sum_{|t|<n}(n-1-|t|)\rho_H(t)^q - \sum_{t\in\Z} \rho_H(t)^q
\right)\\
&=&\frac{1}{q!\sigma^2}\left( -\frac1n\sum_{|t|<n}(|t|+1)\rho_H(t)^q
- \sum_{|t|\geq n} \rho_H(t)^q
\right)\\
&\trianglelefteqslant& \frac1n\sum_{t=1}^{n-1} t^{2qH-2q+1} + \sum_{t=
n}^{\infty} t^{2qH-2q} \trianglelefteqslant n^{-1}+n^{2qH-2q+1}.
\end{eqnarray*}

Now, we assume that $r\leq q-2$ is fixed. We have
\begin{eqnarray*}
E|A_r(n)|^2&=&c(H,r,q)n^{4qH-2}\sum_{i,j,k,\ell=0}^{n-1}
\langle \delta_{k/n},\delta_{\ell/n}\rangle_\HH^{r+1}
\langle\delta_{i/n},\delta_{j/n}\rangle_\HH^{r+1}\\
&&\hskip4cm
\times\langle
\delta_{k/n}^{\otimes q-1-r}\widetilde{\otimes}\delta_{\ell/n}^{q-1-r},
\delta_{i/n}^{\otimes q-1-r}\widetilde{\otimes}\delta_{j/n}^{q-1-r}
\rangle_{\HH^{\otimes 2q-2-2r}}\\
&=&\sum_{\substack{\alpha,\beta\geq 0\\\alpha+\beta=q-r-1}}\,
\sum_{\substack{\gamma,\delta\geq 0\\\gamma+\delta=q-r-1}}
c(H,r,q,\alpha,\beta,\gamma,\delta)
\,B_{r,\alpha,\beta,\gamma,\delta}(n)
\end{eqnarray*}
where $c(\cdot)$ denotes a generic constant depending only on the
objects inside its argument, and
\begin{eqnarray*}
B_{r,\alpha,\beta,\gamma,\delta}(n)&=&
n^{4qH-2}\sum_{i,j,k,\ell=0}^{n-1}
\langle \delta_{k/n},\delta_{\ell/n}\rangle_\HH^{r+1}
\langle \delta_{i/n},\delta_{j/n}\rangle_\HH^{r+1}
\langle \delta_{k/n},\delta_{i/n}\rangle_\HH^{\alpha}\\
&&\hskip5cm\times
\langle \delta_{k/n},\delta_{j/n}\rangle_\HH^{\beta}
\langle \delta_{\ell/n},\delta_{i/n}\rangle_\HH^{\gamma}
\langle \delta_{\ell/n},\delta_{j/n}\rangle_\HH^{\delta}\\
&=& n^{-2}\sum_{i,j,k,\ell=0}^{n-1} \rho_H(k-\ell)^{r+1}
\rho_H(i-j)^{r+1} \rho_H(k-i)^{\alpha}\\
 &&\hskip5cm\times \rho_H(k-j)^{\beta}
\rho_H(\ell-i)^{\gamma} \rho_H(\ell-j)^{\delta}.
\end{eqnarray*}

When $\alpha,\beta,\gamma,\delta$ are fixed, we can decompose the
sum $\sum_{i,j,k,\ell}$ appearing in $B_{r,\alpha,\beta,\gamma,\delta}(n)$
just above,
as follows:
\begin{eqnarray*}
&&\sum_{i=j=k=\ell}
+\left(
\sum_{\substack{i=j=k\\ \ell\neq i}}
+\sum_{\substack{i=j=\ell\\ k\neq i}}
+\sum_{\substack{i=k=\ell\\ j\neq i}}
+\sum_{\substack{j=k=\ell\\ i\neq j}}
\right)
+\left(
\sum_{\substack{i=j, k=\ell\\k\neq i}}
+\sum_{\substack{i=k, j=\ell\\j\neq i}}
+\sum_{\substack{i=\ell, j=k\\j\neq i}}
\right)\\
&&+\left( \sum_{\substack{i=j, k\neq i\\k\neq \ell, \ell\neq i}}
+\sum_{\substack{i=k, j\neq i\\j\neq \ell, k\neq \ell}}
+\sum_{\substack{i=\ell, k\neq i\\ k\neq j, j\neq i}}
+\sum_{\substack{j=k, k\neq i\\k\neq \ell, \ell\neq i}}
+\sum_{\substack{j=\ell, k\neq i\\k\neq \ell, \ell\neq i}}
+\sum_{\substack{k=\ell, k\neq i\\k\neq j, j\neq i}} \right)
+\sum_{\substack{i, j, k, \ell \\ \text{are all different}}}
\end{eqnarray*}
(all these sums must be understood as being defined over indices
$\{i,j,k,\ell\}\!\in\!\{\!0,...,n\!-\!1\!\}^4$). Now, we will deal
with each of these fifteen sums separately.

The first sum is particularly easy to handle: indeed, it is immediately checked that
\begin{eqnarray*}
n^{-2}\sum_{i=j=k=\ell} \rho_H(k-\ell)^{r+1} \rho_H(i-j)^{r+1}
\rho_H(k-i)^{\alpha} \rho_H(k-j)^{\beta} \rho_H(\ell-i)^{\gamma}
\rho_H(\ell-j)^{\delta}\trianglelefteqslant n^{-1}.
\end{eqnarray*}

For the second sum, one can write
\begin{eqnarray*}
&&n^{-2}
\sum_{\substack{i=j=k\\ \ell\neq i}}
\rho_H(k-\ell)^{r+1}
\rho_H(i-j)^{r+1}
\rho_H(k-i)^{\alpha}
\rho_H(k-j)^{\beta}
\rho_H(\ell-i)^{\gamma}
\rho_H(\ell-j)^{\delta}\\
&\trianglelefteqslant& n^{-2}\sum_{i\neq\ell} \rho_H(\ell-i)^{q} \trianglelefteqslant n^{-1}
\sum_{\ell=1}^{n-1} \ell^{2qH-2q} = n^{-1}+n^{2qH-2q}\quad\mbox{by
Lemma \ref{lm-int}}.
\end{eqnarray*}
For the third sum, we can proceed analogously and we again obtain
the bound $n^{-1}+n^{2qH-2q}$.

For the fourth sum, we write
\begin{eqnarray*}
&&n^{-2}
\sum_{\substack{i=k=\ell\\ j\neq i}}
\rho_H(k-\ell)^{r+1}
\rho_H(i-j)^{r+1}
\rho_H(k-i)^{\alpha}
\rho_H(k-j)^{\beta}
\rho_H(\ell-i)^{\gamma}
\rho_H(\ell-j)^{\delta}\\
&\trianglelefteqslant& n^{-2}\sum_{i\neq j} \rho_H(j-i)^{r+1+\beta+\delta}
\trianglelefteqslant n^{-2}\sum_{i\neq j}|j-i|^{(r+1+\beta+\delta)(2H-2)}
\trianglelefteqslant n^{-2}\sum_{i\neq j} |j-i|^{2H-2}\\
&\trianglelefteqslant& n^{-1}\sum_{j=1}^{n-1} j^{2H-2}\trianglelefteqslant n^{-1}+n^{2H-2}
\end{eqnarray*}
(we used the fact that $r+1+\beta+\delta\geq 1$ since
$r,\beta,\delta\geq 0$). For the fifth sum, we can proceed
analogously and we again obtain the bound $n^{-1}+n^{2H-2}$.

For the sixth sum, we have
\begin{eqnarray*}
&&n^{-2}
\sum_{\substack{i=j\\ k=\ell\\k\neq i}}
\rho_H(k-\ell)^{r+1}
\rho_H(i-j)^{r+1}
\rho_H(k-i)^{\alpha}
\rho_H(k-j)^{\beta}
\rho_H(\ell-i)^{\gamma}
\rho_H(\ell-j)^{\delta}\\
&\trianglelefteqslant& n^{-2} \sum_{k\neq i} \rho_H(k-i)^{2q-2-2r} \trianglelefteqslant
n^{-2}\sum_{k\neq i} |k-i|^{(2q-2-2r)(2H-2)}
\trianglelefteqslant n^{-2}\sum_{k\neq i}|k-i|^{4H-4}\\
&\trianglelefteqslant& n^{-1}\sum_{k=1}^{n-1} k^{4H-4} \trianglelefteqslant n^{-1}+n^{4H-4}
\end{eqnarray*}
(here, we used $r\leq q-2$). For the seventh and the eighth sums,
we can proceed analogously and we also obtain $n^{-1}+n^{4H-4}$
for bound.

For the ninth sum, we have
\begin{eqnarray*}
&&n^{-2}
\sum_{\substack{i=j, k\neq i\\k\neq \ell, \ell\neq i}}
\rho_H(k-\ell)^{r+1}
\rho_H(i-j)^{r+1}
\rho_H(k-i)^{\alpha}
\rho_H(k-j)^{\beta}
\rho_H(\ell-i)^{\gamma}
\rho_H(\ell-j)^{\delta}\\
&\trianglelefteqslant&n^{-2} \sum_{\substack{k\neq i\\k\neq \ell, \ell\neq i}}
\rho_H(k-\ell)^{r+1} \rho_H(k-i)^{q-r-1} \rho_H(\ell-i)^{q-r-1}.
\end{eqnarray*}
Now, let us decompose the sum $\sum_{k\neq i, k\neq \ell, \ell\neq i}$ into
$$
\sum_{k>\ell>i}+\sum_{k>i>\ell}+\sum_{\ell>i>k}+\sum_{\ell>k>i}+
\sum_{i>\ell>k}+\sum_{i>k>\ell}.
$$
For the first term (for instance), we have
\begin{eqnarray*}
&&n^{-2}
\sum_{k>\ell>i}
\rho_H(k-\ell)^{r+1}
\rho_H(k-i)^{q-r-1}
\rho_H(\ell-i)^{q-r-1}\\
&\trianglelefteqslant&n^{-2} \sum_{k>\ell>i}
(k-\ell)^{(r+1)(2H-2)}(k-i)^{(q-r-1)(2H-2)}
(\ell-i)^{(q-r-1)(2H-2)}\\
&\trianglelefteqslant&n^{-2} \sum_{k>\ell>i} (k-\ell)^{q(2H-2)}
(\ell-i)^{(q-r-1)(2H-2)}\quad\mbox{since $k-i>k-\ell$}\\
&=&n^{-2}\sum_k \sum_{\ell< k}
(k-\ell)^{q(2H-2)} \sum_{i<\ell}(\ell-i)^{(q-r-1)(2H-2)}\\
&\trianglelefteqslant&n^{-2}\sum_k \sum_{\ell< k}
(k-\ell)^{q(2H-2)} \sum_{i<\ell}(\ell-i)^{2H-2}\quad\mbox{since $q-r-1\ge 1$}\\
&\trianglelefteqslant& n^{-1}\sum_{\ell=1}^{n-1} \ell^{2qH-2q} \sum_{i=1}^{n-1} i^{2H-2}\\
&\trianglelefteqslant&n^{-1}(1+n^{2qH-2q+1})(1+n^{2H-1})\trianglelefteqslant
n^{-1}+n^{2H-2}\quad\mbox{since $2qH-2q+1<0$}.
\end{eqnarray*}
We obtain the same bound for the other terms. By proceeding in the
same way than for the ninth term, we also obtain the bound
$n^{-1}+n^{2H-2}$ for the tenth, eleventh, twelfth, thirteenth and
fourteenth terms.

\noindent For the fifteenth (and last!) sum, we decompose
$\sum_{(i, j, k, \ell\: \text{ are all different})}$ as follows
\begin{equation}\label{anpodipo}
\sum_{k>\ell>i>j}+\sum_{k>\ell>j>i}+\ldots.
\end{equation}
For the first term, we have:
\begin{eqnarray*}
&&n^{-2}
\sum_{k>\ell>i>j}
\rho_H(k-\ell)^{r+1}
\rho_H(i-j)^{r+1}
\rho_H(k-i)^{\alpha}
\rho_H(k-j)^{\beta}
\rho_H(\ell-i)^{\gamma}
\rho_H(\ell-j)^{\delta}\\
&\trianglelefteqslant&n^{-2} \sum_{k>\ell>i>j} (k-\ell)^{q(2H-2)}
(i-j)^{(r+1)(2H-2)}
(\ell-i)^{(q-r-1)(2H-2)}\\
&=&n^{-2}\sum_k \sum_{\ell< k}
(k-\ell)^{q(2H-2)} \sum_{i<\ell}(\ell-i)^{(q-r-1)(2H-2)} \sum_{j<i} (i-j)^{(r+1)(2H-2)}
\\
&\trianglelefteqslant& n^{-1}\sum_{\ell=1}^{n-1}
\ell^{q(2H-2)}\sum_{i=1}^{n-1}i^{(q-r-1)(2H-2)}
\sum_{j=1}^{n-1} j^{(r+1)(2H-2)}\\
&\trianglelefteqslant&n^{-1}(1+n^{2qH-2q+1})(1+n^{(q-r-1)(2H-2)+1})(1+n^{(r+1)(2H-2)+1})\\
&\trianglelefteqslant&n^{-1}(1+n^{2H-1}+n^{2qH-2q+2})\quad\mbox{since $2qH-2q+1<0$ and $r+1,q-r-1\geq 1$}\\
&\trianglelefteqslant&n^{-1}+n^{2H-2}+n^{2qH-2q+1}.
\end{eqnarray*}
The same bound also holds for the other terms in (\ref{anpodipo}).
By combining all these bounds, we obtain
$$
\max_{r=1,\ldots,q-1} E|A_r(n)|^2 \trianglelefteqslant
n^{-1}+n^{2H-2}+n^{2qH-2q+1},
$$
that finally gives:
$$
E\left(\frac1q\|DZ_n\|^2_\HH -1\right)^2 \trianglelefteqslant
n^{-1}+n^{2H-2}+n^{2qH-2q+1}.
$$
The proof of Theorem \ref{BM-rev} is now completed by means of
Proposition \ref{Gauss_Wiener}. \fin


\section{Some remarks about $\chi^2$ approximations}\label{S : Gamma}

The following statement illustrates a natural application of the
results about $\chi^2$ approximations (as discussed in Section
\ref{SS : GammaAPP}) in order to obtain upper bounds in
non-central limit theorems for multiple integrals. Observe that we
focus on double integrals but, at the cost of some heavy notation,
everything can be straightforwardly extended to the case of
integrals of any order $q\geq2$. Recall that the class of functions
$\mathscr{H}_1$ is defined in formula (\ref{H22}).

\begin{prop}
Let $F_n = I_2(f_n)$, $n\geq 1$,
where $f_n \in \HH^{\odot 2}$, be a sequence of double Wiener-Itô
integrals.
Suppose that
$E(F_n^2)\longrightarrow 1$ and $\|f_n \otimes_1
f_n\|_{\HH^{\otimes 2}}\longrightarrow 0$ as $n\rightarrow \infty$.
Then, by defining
$$
H_n = I_4\big(f_n \widetilde{\otimes} f_n\big), \,\, n\geq 1,
$$
one has that
\begin{equation}\label{vabbevabbe}
E\left[\left(2+2H_n - \frac14\|DH_n\|^2_\HH\right)^2\right] \longrightarrow 0,\quad\mbox{as $n\to\infty$,}
\end{equation}
and
\begin{equation}\label{vabbe}
d_{\mathscr{H}_1}(F_n^2-1, N^2-1) \leq 8\sqrt{2} \|f_n \otimes_1
f_n\|_{\HH^{\otimes 2}} + \sqrt{2\pi E\left[\left(2+2H_n -
\frac14\|DH_n\|^2_\HH\right)^2\right]},
\end{equation}
where $N\sim\mathscr{N}(0,1)$.
\end{prop}
{\bf Proof}. First, we have that $F_n\stackrel{\rm Law}{\longrightarrow} N$
by Theorem \ref{david-joe}.
Now,
use the multiplication formula (\ref{multiplication}) to deduce
that
$$
F_n^2-1 = 8\,I_2(f_n \otimes_1 f_n)+H_n.
$$
Since
$$
E\big(I_2(f_n \otimes_1 f_n)^2\big)=2\,\|f_n \otimes_1 f_n\|_{\HH^{\otimes
2}}^2 \longrightarrow 0,\quad\mbox{as $n\to\infty$},
$$
we infer that $H_n\stackrel{\rm Law}{\longrightarrow} N^2-1$, and
therefore that (\ref{vabbevabbe}) must take place, due to Theorem
\ref{ivan-joe}. By the definition of the class $\mathscr{H}_1$,
one also deduces that
$$
d_{\mathscr{H}_1}(F_n^2-1, N^2-1) \leq d_{\mathscr{H}_1}(H_n,
N^2-1) + 8\,E\big|I_2(f_n \otimes_1 f_n)\big|.
$$
The final result is obtained by combining (\ref{K2})--
(\ref{WienerboundG}) with the relation $$E|I_2(f_n \otimes_1
f_n)|\leq \sqrt{E\big(I_2(f_n \otimes_1 f_n)^2\big)}.$$ \fin

We conclude this section with a simple example, showing how one
can apply our techniques to deduce bounds in a non-central limit
theorem, involving quadratic functionals of i.i.d. Gaussian random
variables.

\medskip

\noindent {\bf Example}. Let $(G_k)_{k\geq 0}$ be a sequence of
centered i.i.d. standard Gaussian random variables. Also, let
$(a_k)_{k\in\Z}$ be a sequence of real numbers such that $$a(0)=1,
\quad
a(r)=a(-r),\,r\in\Z,\quad\mbox{and}\quad\sum_{r\in\Z}|a(r)-1|<\infty.$$
In particular, this implies that $a$ is bounded (say, by
$\|a\|_\infty$). Set
$$
F_n=\frac1n\sum_{k,l=0}^{n-1} a(k-l)\big(G_kG_l-\delta_{kl}\big)
,\quad n\geq 1,
$$
where $\delta_{kl}$ denotes the Kronecker symbol. We claim that
$F_n\underset{n\to\infty}{\overset{\rm Law}{\longrightarrow}}
N^2-1$ with $N\sim\mathscr{N}(0,1)$, and our aim is to associate a
bound with this convergence. Observe first that, without loss of
generality, we can assume that $G_k=B_{k+1}-B_k$ where $B$ is a
standard Brownian motion (that $B$ can be therefore regarded as an
isonormal process over $\HH={\rm L}^2(\R_+,dx)$). We then have
$$
DF_n=\frac1n\sum_{k,l=0}^{n-1} a(k-l)\big(G_k {\bf 1}_{[l,l+1]} +
G_l {\bf 1}_{[k,k+1]} \big)
$$
so that
\begin{eqnarray*}
\|DF_n\|^2_{L^2}
&=&\frac1{n^2}\sum_{i,j,k,l=0}^{n-1} a(k-l)a(i-j)\big\langle G_k
{\bf 1}_{[l,l+1]} + G_l {\bf 1}_{[k,k+1]}, G_i {\bf 1}_{[j,j+1]} +
G_j {\bf 1}_{[i,i+1]}
\big\rangle_{L^2}\\
&=&\frac1{n^2}\sum_{i,j,k,l=0}^{n-1} a(k-l)a(i-j)\big(
G_kG_i{\delta}_{lj}+G_kG_j{\delta}_{li}+G_lG_i{\delta}_{kj}+G_lG_j{\delta}_{ik}
\big)\\
&=&\frac4{n}\sum_{k,l=0}^{n-1} G_kG_l\\
&& +\frac1{n^2}\sum_{i,j,k,l=0}^{n-1}
\big(a(k-l)a(i-j)-1\big)\,\big(
G_kG_i{\delta}_{lj}+G_kG_j{\delta}_{li}+G_lG_i{\delta}_{kj}+G_lG_j{\delta}_{ik}
\big).
\end{eqnarray*}
Hence
$$
\frac12\|DF_n\|^2_{L^2}-2F_n-2=A_n+B_n
$$
with
\begin{eqnarray*}
A_n&=&\frac2{n}\sum_{k,l=0}^{n-1}\big(1-a(k-l)\big)
G_kG_l\\
B_n&=&\frac1{2n^2}\sum_{i,j,k,l=0}^{n-1}
\big(a(k-l)a(i-j)-1\big)\,\big(
G_kG_i{\delta}_{lj}+G_kG_j{\delta}_{li}+G_lG_i{\delta}_{kj}+G_lG_j{\delta}_{ik}
\big).
\end{eqnarray*}
We have
\begin{eqnarray*}
E(A_n^2)&=&\frac4{n^2}\sum_{i,j,k,l=0}^{n-1}\big(1-a(k-l)\big)\big(1-a(i-j)\big)
E\big(G_kG_lG_iG_j\big)\\
&=&\frac8{n^2}\sum_{i,k=0}^{n-1}\big(1-a(k-i)\big)^2
\leq\frac8n\big(1+\|a\|_\infty\big)\sum_{r\in\Z}\big|1-a(r)\big|=O(1/n).
\end{eqnarray*}
On the other hand, we have
$$
B_n=B_n^1+B_n^2+B_n^3+B_n^4
$$
with
\begin{eqnarray*}
B^1_n&=&\frac1{2n^2}\sum_{i,j,k,l=0}^{n-1}
\big(a(k-l)a(i-j)-1\big)\,G_kG_i{\delta}_{lj}
\end{eqnarray*}
and similar computations hold for the other terms. Observe that
\begin{eqnarray*}
B^1_n&=&\frac1{2n^2}\sum_{i,j,k=0}^{n-1}
\big(a(k-j)a(i-j)-1\big)\,G_kG_i\\
&=&\frac1{2n^2}\sum_{i,k=0}^{n-1}\alpha_{ki}\,G_kG_i,\quad\mbox{with
$\alpha_{ki}=\sum_{j=0}^{n-1}\big(a(k-j)a(i-j)-1\big)$}.
\end{eqnarray*}
We have
\begin{eqnarray*}
\big|\alpha_{ki}\big|=\left|\sum_{j=0}^{n-1}\big(a(k-j)-1\big) +
\sum_{j=0}^{n-1}a(k-j)\big(a(i-j)-1\big)\right|\leq
(1+\|a\|_\infty)\sum_{r\in\Z}|a(r)-1|.
\end{eqnarray*}
Consequently
$$
E\big|B^1_n\big|^2=\frac1{4n^4}\sum_{i,j,k,l=0}^{n-1}\alpha_{ki}\alpha_{jl}\,E\big(G_kG_iG_lG_j\big)=O(1/n^2).
$$
Similarly, the same bound holds for $E\big|B^i_n\big|^2$,
$i=2,3,4$. Finally, by combining all the previous estimates, we
obtain
$$
\sqrt{E\big(\frac12\|DF_n\|^2_{L^2}-2F_n-2\big)^2}=O(1/\sqrt{n}),
$$
and therefore, by using Theorem \ref{gam-thm} and the fact that
$N^2-1 \stackrel{\rm Law}{=} F(1)$, we deduce that there exists a
positive constant $C>0$ (independent of $n$) such that
$$
d_{\mathscr{H}_1}(G,N^2 - 1)\leq C /\sqrt{n},
$$
where the class $\mathscr{H}_1$ is defined in (\ref{H22}).
\fin

\section{An attempt at unification}\label{S : Unif}
In this section, we show that the computations contained in
Section \ref{SS : GaussAPP} and Section \ref{SS : GammaAPP},
respectively, in the Gaussian case and the Gamma case, can be
unified, by means of the general theory of approximations
developed by Ch. Stein in \cite[Lecture VI]{Stein_book}.

Let $Z$ be a real-valued random variable having an absolutely
continuous distribution with density $p(x)$, $x\in \R$. We make
the following assumptions:
\begin{enumerate}
\item[\textbf{(A1)}] $Z$ is integrable and centered, that is,
\begin{equation}\label{centered}
E|Z|<\infty \quad \text{and} \quad
E(Z)=\int^{+\infty}_{-\infty}yp(y)dy=0;
\end{equation}
\item[\textbf{(A2)}] there exist (possibly infinite) numbers $a,b$ such that
$-\infty \leq a <0<b\leq +\infty$, and the support of the density
$p$ exactly coincides with the open interval $(a,b)$, that is,
\begin{equation}\label{support}
p(x)>0 \quad \text{if, and only if,}\quad x\in(a,b).
\end{equation}
\end{enumerate}

\begin{rem}
{\rm At the cost of some heavier notation, one could easily
generalize the results of this section, in order to accommodate
the case of a density $p$ whose support is a union of open (and
possibly infinite) intervals.}
\end{rem}

With a $Z$ verifying assumptions \textbf{(A1)-(A2)}, we associate
the real-valued mapping $\tau(\cdot)$, defined as

\begin{equation}\label{tau}
x \mapsto \tau(x)= \frac{\int^\infty_x
yp(y)dy}{p(x)}\mathbf{1}_{x\in(a,b)} = -\frac{\int_{-\infty}^x
yp(y)dy}{p(x)}\mathbf{1}_{x\in(a,b)}, \quad x\in \R.
\end{equation}

Note that $\tau$ is well-defined on $\R$, due to assumptions
(\ref{centered})-(\ref{support}). Also, relation (\ref{centered})
implies that $\tau(x) \geq 0$ for every $x$ and $\tau(x) > 0$ if,
and only if, $x\in (a,b)$. The following result, which is proved
in \cite{Stein_book}, states that, under some additional
assumptions, the mapping $\tau$ completely characterizes the
density $p$, and therefore the law of $Z$.

\begin{lemme}[Lemma 3, p. 61 in \cite{Stein_book}]\label{L : TAU}
Let the reals $a,b$ be such that $-\infty\leq a <0 < b \leq
+\infty$, and consider a continuous function $\tau(\cdot)\geq0$ on
$\R$ such that
\begin{equation}\label{tauSupp}
\tau(x) > 0 \quad \text{if, and only if,} \quad x\in (a,b).
\end{equation}
Then, if
\begin{equation}\label{tauexplosion}
\int_0^b(y/\tau(y))dy=+\infty \quad \text{and} \quad \int_a^0
(y/\tau(y))dy=-\infty,
\end{equation}
there exists a unique (up to sets of zero Lebesgue measure)
probability density $p_\tau(\cdot)$ on $\R$ such that the support
of $p_\tau$ exactly coincides with the interval $(a,b)$ and
\begin{equation} \label{tauCaract}
\int^{+\infty}_{-\infty}yp_\tau(y)dy=0 \quad \text{and} \quad
\tau(x)=
\frac{\int^{+\infty}_{x}yp_\tau(y)dy}{p_\tau(x)}\mathbf{1}_{x\in(a,b)},
\quad x\in\R.
\end{equation}
The explicit form of $p$ is given by
\begin{equation}\label{taudensity}
p_\tau(x) = \frac{1}{C}\times
\frac{\mathrm{e}^{-\int_0^x\frac{ydy}{\tau(y)}}}{\tau(x)}\mathbf{1}_{x\in(a,b)},
\quad x\in \R,
\end{equation}
where $ C = \int^b_a
\frac{\mathrm{e}^{-\int_0^x\frac{ydy}{\tau(y)}}}{\tau(x)}dx, $ and
we used the notational convention $\int_0^x = - \int_x^0$ whenever
$x<0$.
\end{lemme}

We will see later on that property (\ref{tauexplosion}) is
verified by the functions $\tau$ associated with densities in the
Pearson's family of continuous distributions. Now let $X$ have a
density $p$ verifying assumptions (\textbf{A1})-(\textbf{A2})
above, and let $\tau$ be the mapping given by (\ref{tau}) (for the
time being, we do not suppose that (\ref{tauexplosion}) is
verified). We define the \textsl{Stein operator} $T_\tau$,
associated with $p$ and $\tau$, as the differential operator
\begin{equation}\label{SteinOp}
T_\tau f(x) = \tau(x)f'(x)-xf(x), \quad x\in\R,
\end{equation}
acting on differentiable functions $f$. Now fix a function $h$ which is piecewise
continuous on $\R$ and such that $E(h(Z))$ is well-defined. The
\textsl{Stein equation}, associated with $p$, $\tau$ and $h$, is
the first order differential equation
\begin{equation}\label{SteinEQ}
h(x)-E(h(Z)) = T_\tau f(x), \quad x\in\R,
\end{equation}
where $T_\tau f$ is defined in (\ref{SteinOp}). If $\tau$ verifies
(\ref{tauexplosion}), then (due to Lemma \ref{L : TAU})
$E(Z)=E_\tau (h)$, where $E_\tau (h) = \int h(y)p_\tau(y)dy,$ and
$p_\tau$ is the density given by (\ref{taudensity}). It follows
that, in this case, one can rewrite (\ref{SteinEQ}) as
\begin{equation}\label{SteinEQ2}
h(x)-E_\tau (h) = T_\tau f(x), \quad x\in\R,
\end{equation}
in order to emphasize the role of $\tau$. The next result, whose
(rather straightforward) proof is once again given by Stein in
\cite{Stein_book}, states that, under (\ref{tauexplosion}), the
equation (\ref{SteinEQ}) admits a unique continuous and bounded
solution.

\begin{lemme}[Lemma 4, p. 62 in \cite{Stein_book}]\label{L : Stein
sol} Let $\tau$ satisfy (\ref{tauSupp}) and (\ref{tauexplosion}),
and let $p_\tau$ be the density associated with $\tau$
\textit{via} (\ref{taudensity}). Then, since $\tau$ has support in
$(a,b)$, every solution $f$ of (\ref{SteinEQ2}) must necessarily
be such that
\begin{equation}\label{SteinTails}
f(x)=\frac{h(x)-E_\tau (h)}{x}, \quad x\in \R  \backslash (a,b).
\end{equation}
Moreover, whenever $h$ is bounded and piecewise continuous, the
equation (\ref{SteinEQ2}) admits a unique solution $f$ which is
bounded and continuous on $(a,b)$. This unique solution is defined
by (\ref{SteinTails}) on $\R  \backslash (a,b)$, and by
\begin{equation}\label{SteinSOLUTION}
f(x)=\int^x_a (h(y)-E_\tau
(h))\frac{\mathrm{e}^{\int^x_y\frac{zdz}{\tau(z)}}}{\tau(y)}dy,
\quad \text{for every} \,\,\,\, x\in (a,b).
\end{equation}
\end{lemme}
Given a bounded and piecewise continuous function $h$ on $\R$, we
define the function $U_\tau h$ as
\begin{eqnarray}\label{SteinU}
U_\tau h(x)& =&\left\{
\begin{array}{ll}
\frac{h(x)-E_\tau (h)}{x}, & \text{if } x\in \R  \backslash (a,b) \\
\int^x_a (h(y)-E_\tau
(h))\frac{\mathrm{e}^{\int^w_y\frac{zdz}{\tau(z)}}}{\tau(y)}dy, &
\text{if }x\in (a,b),
\end{array}
\right.
\end{eqnarray}
so that one can rephrase Lemma \ref{L : Stein sol} by saying that
$U_\tau h$ is the unique solution of (\ref{SteinEQ2}) which is
bounded and continuous on $(a,b)$ (note that $U_\tau h$ can be
discontinuous only at $a$ or $b$, whenever they are finite). We
also record the following consequence of the calculations
contained in \cite[formulae (34)-(35), pp. 64-65]{Stein_book}:
\textsl{if $h$ is bounded and piecewise continuous, then}
\begin{equation}\label{SteinBOUND}
\sup_{x\in (a,b)}[|x U_\tau h(x)|+|\tau(x) U'_\tau h(x)|]\leq
6\sup_{x\in (a,b)}|h(x)|,
\end{equation}
where $U'_\tau h = (U_\tau h)'$. Note that, due to
(\ref{SteinTails}), one deduces immediately from
(\ref{SteinBOUND}) that
\begin{equation}\label{SteinBOUND2}
\sup_{x\in \R}[|x U_\tau h(x)|+|\tau(x) U'_\tau h(x)|]\leq K
\sup_{x\in \R}|h(x)|,
\end{equation}
where $K = 2\max\{3;1/a;1/b\}$ (with $1/\pm\infty =0$). The next
statement provides a typical `Stein-type characterization' of the
law of $Z$. It is a general version of Lemma
\ref{Stein_Lemma_Gauss}(i) and Lemma \ref{Stein_Lemma_Gamma}(i).

\begin{prop}\label{P : STeinCHAR}
Let $Z$ be a random variable having a density $p$ verifying
assumptions (\textbf{A1}) and (\textbf{A2}). Let $\tau$ be related
to $p$ by (\ref{tau}).
\begin{enumerate}
\item[(i)] For every differentiable $f$ such that
$E|\tau(Z)f'(Z)|<\infty$, one has that $E|Zf(Z)|<\infty$ and
\begin{equation}\label{SteinCHAR}
E[T_\tau(Z)] = E[\tau(Z)f'(Z)-Zf(Z)] = 0.
\end{equation}
\item[(ii)]Suppose in addition that $\tau$ verifies (\ref{tauexplosion}).
Let $Y$ be a real-valued random variable with an absolutely
continuous distribution. Suppose that, for every differentiable
$f$ such that the mapping $x \mapsto |\tau(x)f'(x)| + |xf(x)|$
($x\in\R$) is bounded, one has that
\begin{equation}\label{SteinChar2}
E[T_\tau(Y)]=E[\tau(Y)f'(Y)-Yf(Y)] = 0.
\end{equation}
Then, $Y\stackrel{\rm Law}{=} Z$.
\end{enumerate}
\end{prop}
\textbf{Proof. } Part (i) is proved in \cite[Lemma 1, p.
69]{Stein_book}. Part (ii) is a consequence of the fact that, if
(\ref{SteinChar2}) is in order, then (due to
(\ref{SteinSOLUTION})--(\ref{SteinBOUND})), for every bounded and
piecewise continuous function $h$ on $\R$, $0 = E[\tau(Y)U'_\tau
h(Y)\!-\!YU_\tau h(Y)]\! =\!E[h(Y)]\!-\!E_\tau
(h)\!=\!E[h(Y)]\!-\!E[h(Z)].$ \fin

The following corollary can be proved along the lines of Theorem
\ref{main-tool} and Theorem \ref{gam-thm}.

\begin{cor}
Let $Z$ be a random variable having a density $p$ verifying
assumptions (\textbf{A1}) and (\textbf{A2}). Let $\tau$ be related
to $p$ by (\ref{tau}). Let $F\in\sk^{1,2}$ be a smooth functional
of some isonormal Gaussian process. Assume moreover that $E(F)=0$ and the law
of $F$ is absolutely continuous with respect to the Lebesgue
measure. Then, for every bounded and piecewise continuous function
$h$, we have
\begin{eqnarray}
E(h(F))-E(h(Z))& =& E[\tau(F)(U_\tau h)'(F)-F U_\tau h(F)] \label{GenEQ}\\
&=&E[(U_\tau h)'(F)(\tau(F)-\langle DF,-DL^{-1}F\rangle_\HH)]
\label{GenEQ2}.
\end{eqnarray}
Also,
\begin{equation}\label{GenINEQ}
|E(h(F))-E(h(Z))|\leq E[(U_\tau h)'(F)^2]^{1/2}E[(\tau(F)-\langle
DF,-DL^{-1}F\rangle_\HH)^2]^{1/2}.
\end{equation}
\end{cor}

It is not difficult to see that the conclusions of Theorem
\ref{main-tool} and Theorem \ref{gam-thm} are indeed corollaries
of formula (\ref{GenINEQ}), corresponding, respectively, to
$\tau(x) = 1$ and $\tau(x) = 2(x+\nu)_+$. Plainly, a study of
general expressions such as the RHS of (\ref{GenINEQ}) would
require a fine analysis of the properties of the solutions to the
Stein equation (\ref{SteinEQ}) (similar to the ones performed in
the Gamma case by Luk and Pickett, respectively, in \cite{Luk} and
\cite{These_Pickett}). This topic is clearly outside the scope of
the present paper. However, we conjecture that such a study could
be successfully performed in the case where the density $p$
belongs to the Pearson's family of curves. Indeed, in this case
the function $\tau$ can be neatly characterized in terms of
polynomials of degree 2.

To see this, let $Z$ satisfy \textbf{(A1)-(A2)}, and let $\tau$
satisfy (\ref{tauexplosion}). We say that $Z$ is a (centered)
member of the \textsl{Pearson's family} of continuous
distributions, whenever the density $p=p_\tau$ (see
(\ref{taudensity})) satisfies the differential equation
\begin{equation}\label{Pearson}
\frac{p'(x)}{p(x)}=\frac{a_0+a_1 x}{b_0+b_1 x+b_2 x^2}, \quad
x\in(a,b),
\end{equation}
for some real numbers $a_0,a_1,b_0,b_1,b_2$. We refer the reader
e.g. to \cite[Sec. 5.1]{DiacZab} for an introduction to the
Pearson's family. Here, we shall only observe that there are
basically five families of distributions satisfying
(\ref{Pearson}): the centered normal distributions, centered Gamma
and beta distributions, and distributions that are obtained by
centering densities of the type $p(x)=Cx^{-\alpha}\exp(-\beta/x)$
or $p(x)=C(1+x)^{-\alpha}\exp(\beta\arctan(x))$ ($C$ being a
suitable normalizing constant). The next result, proved in
\cite[Theorem 1, p. 65]{Stein_book}, states that a density belongs
to the class of the Pearson's curves if, and only if, its
associated mapping $\tau$ is a polynomial of degree $\leq2$. The
reader is also referred to \cite[Sec. 2 and Sec. 4]{Schoutens} for
several related results and explicit computations involving
orthogonal polynomials.

\begin{thm}[Stein]\label{T : St + Pears} Let $Z$ satisfy \textbf{(A1)-(A2)}, and let $\tau$ satisfy
(\ref{tauexplosion}). Then, the density $p=p_\tau$ is such that
$\tau(x)=\alpha x^2 +\beta x + \gamma$, $x\in(a,b)$ (with
$\alpha,\,\beta, \,\gamma$ constants) if, and only if, $p$
satisfies (\ref{Pearson}) for every $x\in(a,b)$ and for $a_0=\beta
$, $a_1=2\alpha + 1$, $b_0=\gamma$, $b_1=\beta$ and $b_2=\alpha$.
\end{thm}

Of course, in order for (\ref{tauexplosion}) to be satisfied, one
must have that the $\tau(a)=0$ (whenever $a$ is finite) and
$\tau(b)=0$ (whenever $b$ is finite). As already discussed, the
centered Gaussian distribution is a member of the Pearson's
family, corresponding to the case $a = -\infty$, $b=+\infty$ and
$\tau(x) = 1$. Analogously, a centered Gamma random variable
$F(\nu)$ as in (\ref{GammaCent}) has a density of the Pearson
type, with characteristics $a=-\nu$, $b=+\infty$ and $\tau(x) =
2(x+\nu)_+$.

\section{Two proofs}\label{S : Proofs}

\subsection{Proof of Lemma \ref{Stein_Lemma_Gamma}}

\noindent\underline{Proof of Point (i)}. One could use directly
Proposition \ref{P : STeinCHAR} in the case $\tau(x)=2(x+\nu)_+$.
Alternatively, observe first that, for every $\nu>0$, the random
variable $F^*(\nu):=F(\nu)+\nu$ has a non-centered Gamma law with
parameter $\nu/2$. The fact that
$$
E[2F^*(\nu)f'(F^*(\nu)-\nu)]=E[2(F^*(\nu)-\nu)f'(F^*(\nu))],
$$
for every $f$ as in the statement, is therefore an immediate
consequence of \cite[Proposition 1 and Section 4(2)]{Schoutens}.
Now suppose that $W$ verifies (\ref{SteinTruc2}). By choosing $f$
with support in $(-\infty, -\nu)$, one deduces immediately that
$P(W\leq -\nu)=0$. To conclude, we apply once again the results
contained in \cite{Schoutens}, to infer that the relations
$$
P(W\leq -\nu)=0 \quad \text{and} \quad E[2(W+\nu)f'(W)-Wf(W)]=0
$$
imply that, necessarily, $W+\nu \stackrel{\rm Law}{=} F^*(\nu)$.

\noindent\underline{Proof of Point (ii)}. Fix $\nu>0$, consider a
function $h$ as in the statement and define $h_\nu(y)=h(y-\nu)$,
$y>0$. Plainly, $h_\nu$ is twice differentiable, and
$|h_\nu(y)|\leq c \exp\{-\nu a\} \exp\{a y\}$, $y>0$ (recall that
$a>1/2$). In view of these properties, according to Luk \cite[Th.
1]{Luk}, the second-order Stein equation
\begin{equation}\label{StNONCGamma}
h_\nu(y) -E(h_\nu(F^*(\nu))= 2yg''(y)-(y-\nu)g'(y), \quad y>0,
\end{equation}
(where, as before, we set $F^*(\nu)=F(\nu)+\nu$) admits a solution
$g$ such that $\|g'\|_\infty \leq 2 \|h'\|_\infty$ and
$\|g''\|_\infty \leq \|h''\|_\infty$. Since $f(x)=g'(x+\nu)$,
$x>-\nu$, is a solution of (\ref{SteinGammaEq}), the conclusion is
immediately obtained.

\noindent\underline{Proof of Point (iii)}. According to a result
of Pickett \cite{These_Pickett}, as reported in \cite[Lemma
3.1]{Reinert_sur}, when $\nu\geq 1$ is an integer, the ancillary
Stein equation (\ref{StNONCGamma}) admits a solution $g$ such that
$\|g'\|_\infty \leq \sqrt{2\pi/\nu} \|h\|_\infty$ and
$\|g''\|_\infty \leq \sqrt{2\pi/\nu} \|h'\|_\infty$. The
conclusion is obtained as in the proof of Point (ii).
\subsection{Proof of Theorem \ref{quidechire}}
We begin with a technical lemma.
\begin{lemme}
Let $F=I_2(f)$ be a random variable living in the second Wiener
chaos of an isonormal Gaussian process $X$ (over a real Hilbert
space $\HH$). Then
\begin{equation}\label{df4}
E\big(\|DF\|^4_\HH\big)=\frac23E(F^4)+2E(F^2)^2.
\end{equation}
\end{lemme}
{\bf Proof}. Without loss of generality, we can assume that
$\HH={\rm L}^2(A,\mathscr{A},\mu)$, where $(A,\mathscr{A})$ is a
measurable space, and $\mu$ is a $\sigma$-finite and non-atomic
measure. On one hand, thanks to the multiplication formula
(\ref{multiplication}), we can write
$$
F^2=I_4(f\otimes f)+4\,I_2(f\otimes_1 f) + E\big(F^2\big).
$$
In particular, this yields
$$
L(F^2)=-4\, I_4(f\otimes f)-8\,I_2(f\otimes_1 f).
$$
On the other hand, (\ref{dtf}) implies that $ D_a F = 2\,
I_1\big(f(\cdot,a)\big). $ Consequently, again by (\ref{multiplication}):
\begin{eqnarray}
\|DF\|^2_\HH&=&4\int_A I_1\big(f(\cdot,a)\big)^2\mu(da)\nonumber\\
&=&4\int_A I_2\big(f(\cdot,a)\otimes f(\cdot,a)\big)\mu(da)
+ E\big(\|DF\|_\HH^2\big)\nonumber\\
&=&4\,I_2(f\otimes_1 f) + 2 E(F^2),\quad \mbox{by (\ref{MallMoments})
and since $\int_A f(\cdot,a)\otimes f(\cdot,a)\mu(da)=f\otimes_1f$}.\nonumber\\
\label{dfk}
\end{eqnarray}
Taking into account the orthogonality between multiple stochastic
integrals of different orders, we deduce
\begin{equation}\label{sde}
E\big[\|DF\|^2_\HH\,L(F^2)\big]=-32\,E\left[\big(I_2(f\otimes_1
f)\big)^2\right] =-2\,E\left[\|DF\|^2_\HH \big(F^2
-E(F^2)\big)\right].
\end{equation}
Finally, we have
\begin{eqnarray*}
E\big[\|DF\|^4_\HH\big]&=&E\big[\|DF\|^2_\HH\langle DF,DF\rangle_\HH\big]\\
&=&E\big[\|DF\|^2_\HH\big(\delta DF\times F-\frac12\delta D(F^2)\big)\big]\quad\mbox{by identity (\ref{dfu})},\\
&=&2\,E\big[\|DF\|^2_\HH F^2\big]+\frac12\,E\big[\|DF\|^2_\HH\,L(F^2)\big]\quad\mbox{using $\delta D=-L$},\\
&=&E\big[\|DF\|^2_\HH F^2\big]+E(F^2)E\big[\|DF\|^2_\HH\big]\quad\mbox{using (\ref{sde})},\\
&=&\frac23 \,E\big(F^4\big)+2\,E\big(F^2\big)^2\quad\mbox{by (\ref{MallMoments})}.
\end{eqnarray*}
\fin Now, let us go back to the proof of the first point in
Theorem \ref{quidechire}. In view of Theorem \ref{main-tool}, it
is sufficient to prove that
\begin{equation}\label{int}
E\left(\left|1-\frac12\|DZ_n\|^2_\HH\right|^2\right) \leq
\frac16\big|E(Z_n^4)-3\big|+\frac{3+E(Z_n^2)}{2}\,\big|E(Z_n^2)-1\big|.
\end{equation}
We have
\begin{eqnarray*}
&&E\left(\left|1-\frac12\|DZ_n\|^2_\HH\right|^2\right)\\
&&=1-E(\|DZ_n\|^2_\HH)+\frac14E(\|DZ_n\|^4_\HH)\\
&&=1-2E(Z_n^2)+\frac16E(Z_n^4)+\frac12E(Z_n^2)^2\quad\quad\quad\quad\mbox{by (\ref{MallMoments}) and (\ref{df4})}\\
&&=\frac16(E(Z_n^4)-3)+(E(Z_n^2)-1)\left(\frac12E(Z_n^2)-\frac32\right).
\end{eqnarray*}
The estimate (\ref{int}) follows immediately.

\smallskip

\noindent Similarly, for the second point of Theorem
\ref{quidechire}, it is sufficient to prove (see Proposition
\ref{Gamma_Wiener}) that
\begin{eqnarray}\label{int2}
&&E\left(\left|2Z_n-2\nu-\frac12\|DZ_n\|^2_\HH\right|^2\right) \\
&&\hskip2.9 cm\leq
\frac16\big|E(Z_n^4)-12E(Z_n^3)-12\nu^2+48\nu\big|
\!+\!\frac{\big|8-6\nu+E(Z_n^2)\big|}{2}\,\big|E(Z_n^2)-2\nu\big|.\nonumber
\end{eqnarray}
By using the relations
\begin{eqnarray*}
&&E\left(\left|2Z_n-2\nu-\frac12\|DZ_n\|^2_\HH\right|^2\right) \\
&=&4E(Z_n^2)+4\nu^2+\frac14E(\|DZ_n\|^4_\HH)-2E(Z_n\|DZ_n\|^2_\HH)-2\nu E(\|DZ_n\|_\HH^2)\\
&=&4(1-\nu)E(Z_n^2)+4\nu^2+\frac16E(Z_n^4)+\frac12E(Z_n^2)^2-2E(Z_n^3)\quad\mbox{by (\ref{MallMoments}) and (\ref{df4})}\\
&=&(E(Z_n^2)-2\nu)\big(4-3\nu+\frac12E(Z_n^2)\big)+\frac16\big(
E(Z_n^4)-12E(Z_n^3)-12\nu^2+48\nu \big),
\end{eqnarray*}
the estimate (\ref{int2}) follows immediately. \fin

\bigskip

\noindent {\bf Acknowledgments}. We thank an anonymous referee for
interesting suggestions and remarks. We are grateful to G. Reinert
for bringing to our attention references \cite{Luk} and
\cite{These_Pickett}. We also thank J. Dedecker for useful
discussions.

\noindent \textsl{This paper is dedicated to the memory of Livio
Zerbini.}

\end{document}